\documentstyle[12pt]{article}
\title{\LARGE\bf Coloring ordinals by reals}
\date{February 24,\ 1997}
\author{Saka\'e Fuchino%
\ifcommented
\bigskip\bigskip\\
\fbox{\parbox{8.8cm}{\footnotesize\tt\mbox{}\hfill
This is yet a preliminary version of the paper.\hfill\mbox{}\\
\mbox{}\hfill Any comments are appreciated.\hfill\mbox{}
}}
\fi
}
\newif\iftesting
\newif\ifcommented 
\commentedtrue
\iftesting
\let\Label\label%
\def\label#1{\mbox{}\marginpar{{\tiny #1}}\Label{#1}\ignorespaces}%
\fi

\ifcommented
\fi
{\ifcommented\end{footnotesize}\medskip\\\fi}
\newcommand{\memo}[1]{\ifcommented #1\fi}

\newtheorem{Thm}{{\bf Theorem}}[section]
\newtheorem{Cor}[Thm]{{\bf Corollary}}

\newtheorem{Lemma}[Thm]{{\bf Lemma}}

\newtheorem{Claim}{{\bf Claim}}[Thm]
\newtheorem{Problem}[Thm]{{\bf Problem}}

\newcommand{\Thmof}[1]{{Theorem \ref{#1}}}

\newcommand{\Lemmaof}[1]{{Lemma \ref{#1}}}

\newcommand{\Thmabove}{{Theorem \number\theThm}}
\newcommand{\Corabove}{{Corollary \number\theThm}}

\newcommand{\prf}{{\bf Proof\ \ }\ignorespaces}
\newcommand{\prfof}[1]{{\bf Proof of #1\,:\,\ }\ignorespaces}
\newcommand{\prfofClaim}{\raisebox{-.4ex}{\Large $\vdash$\ \ }}
\newsavebox{\qedbox}\sbox{\qedbox}{% QED symbol (by U.Fuchs modified by S.F.)
{\unitlength=0.07mm \begin{picture}(40,60)
\put(0,0){\framebox(30,44)[cc]{}}
\put(30,-7){\rule{7\unitlength}{44\unitlength}}
\put(10,-7){\rule{27\unitlength}{7\unitlength}}
\end{picture}}}
\newcommand{\qed}{\mbox{}\hfill\usebox{\qedbox}}
\newcommand{\smallqed}%
{\mbox{}\smallskip\hfill\raisebox{-.4ex}{\Large $\dashv$}\\}
\newcommand{\qedof}[1]%
{\mbox{} \hspace*{\fill}{\usebox{\qedbox}{~(#1)}}%
\mbox{}}
\newcommand{\Qedof}[1]%
{\mbox{} \hspace*{\fill}{\usebox{\qedbox}%
{~(#1~\number\theThm)}}}
\newcommand{\qedofThm}{\Qedof{Theorem}}
\newcommand{\qedofCor}{\Qedof{Corollary}}

\newcommand{\qedofLemma}{\Qedof{Lemma}}

\newcommand{\qedskip}{\medskip}
\newcommand{\qedofClaim}%
{\mbox{}\hfill\raisebox{-.4ex}{\Large $\dashv$ }\nolinebreak%
\mbox{~(Claim~\number\theClaim)}}

\newcommand{\st}{such that}
\newcommand{\wrt}{with respect to}

\newcommand{\wolog}{without loss of generality}
\newcommand{\Wolog}{Without loss of generality}
\newcommand{\po}{partial ordering}
\newcommand{\pos}{partial orderings}
\newcommand{\tenten}{\ifmmode,\ldots,\else\mbox{\rm,\ldots,} \fi}
\newcommand{\Equivto}{$\,\Leftrightarrow\,$}
\newcommand{\equivto}{\ \Leftrightarrow\ }
\newenvironment{assertion}[1]{\begin{trivlist}
\newbox\assertbox
\dimen255=\textwidth
\setbox\assertbox=\hbox{\hspace*{\parindent}{#1}\hspace{\labelsep}}
\advance\dimen255 by -1\wd\assertbox
\advance\dimen255 by -1ex
\item[]\unhbox\assertbox\hfill
\begin{minipage}[t]{\dimen255}}%
{\end{minipage}\end{trivlist}}
{\hfill\mbox{}\end{trivlist}}
\def\assert#1{\makebox[5ex][l]{\rm (#1)}\ignorespaces}
\def\lassert#1{\llap{\makebox[5ex][l]{\rm (#1)}}\ignorespaces}
\def\assertof#1{{\rm (#1)}}
\newcommand{\assertskip}{\smallskip}
\newenvironment{markedformula}[1]%
	{\begin{trivlist}\item[\,\,#1]\mbox{}\hfill$\displaystyle}%
	{$\hfill\mbox{}\end{trivlist}}
\newcommand{\setof}[2]{\{#1\,:\,#2\}}
\newcommand{\ssetof}[1]{\{#1\}}
\newcommand{\seqof}[2]{\langle#1\,:\,#2\rangle}
\newcommand{\pairof}[1]{\langle#1\rangle}
\newcommand{\psetof}[1]{{\cal P}\/(#1)}
\newcommand{\cardof}[1]{\mathopen{|\,}#1\mathclose{\,|}}
\newcommand{\mapping}[3]{#1:#2\rightarrow #3}
\newcommand{\fnsp}[2]{\mbox{}^{#1^{\mbox{}\!}}#2}
\newcommand{\dotbigcup}{\mathop{\dot{\bigcup}}}
\newcommand{\restr}%
{{\hspace{0.1ex}|\hspace{-0.02ex}{\grave{}}\hspace{0.8ex}}}
\newcommand{\reals}{\mbox{$\rm I_{\!\!}R$}}
\newcommand{\bbbone}{{\mathchoice {\rm 1\mskip-4mu l} {\rm 1\mskip-4mu l}
{\rm 1\mskip-4.5mu l} {\rm 1\mskip-5mu l}}}
\newcommand{\regembed}%
	{\mathrel{{<}\llap{\raisebox{0.2ex}{$\scriptstyle\,\circ$}}}}
\newcommand{\forces}[2]{\,\|\hspace{-.35ex}\mbox{\sf--}_{\,#1\,}%
\mbox{\rm``}\,#2\,\mbox{\rm''}}
\newcommand{\notforces}[2]{\rlap{\ /}\|\hspace{-.35ex}\mbox{\sf--}_{\,#1\,}%
\mbox{\rm``}\,#2\,\mbox{\rm''}}
\newcommand{\lessnoneq}%
{\mathrel{\raisebox{-0.8ex}{$\stackrel{<}{\scriptstyle\,\not=\,}$}}}

\newcommand{\xmbox}[1]{ $\relax{\rm #1}\relax$ }

\newcommand{\dom}{\mathop{\rm dom}}
\newcommand{\range}{\mathop{\rm rng}}
\newcommand{\supp}{\mathop{\rm supp}}
\newcommand{\cof}{\mathop{\rm cof}}
\newcommand{\Fn}{\mathop{\rm Fn}}

\newcommand{\ZFC}{{\rm ZFC}}
\newcommand{\CH}{{\rm CH}}
\newcommand{\GCH}{{\rm GCH}}
\newcommand{\MA}{{\rm MA}}
\newcommand{\HP}{{\rm HP}}
\newcommand{\IP}{{\rm IP}}
\newcommand{\OCA}{{\rm OCA}}
\newcommand{\On}{{\rm On}}

\newcommand{\calH}{{\cal H}}

\newcommand{\dota}{{\dot a}}
\newcommand{\dotf}{{\dot f}}
\newcommand{\dotg}{{\dot g}}
\newcommand{\dotx}{{\dot x}}
\newcommand{\doty}{{\dot y}}

\newcommand{\dotA}{{\dot A}}
\newcommand{\dotC}{{\dot C}}
\newcommand{\dotD}{{\dot D}}
\newcommand{\dotG}{{\dot G}}
\newcommand{\dotS}{{\dot S}}
\newcommand{\dotX}{{\dot X}}

\begin{document}
\maketitle
\begin{abstract}
We introduce several new set-theoretic axioms 
formulated in terms of coloring of ordinals by reals. We show 
that these axioms generalize the axioms considered by I.\ Juh\'asz, L.\ 
Soukup and Z.\ Szentmikl\'ossy \cite{soukup-etal}, and 
give a class of \pos\ including Cohen p.o.-sets which force these axioms. 
\end{abstract}
\section{Introduction}\label{S1}
In \cite{soukup-etal}, I.\ Juh\'asz, L.\ Soukup and Z.\ Szentmikl\'ossy 
introduced several new set-theoretic axioms which hold in Cohen 
models, i.e.\ the models of \ZFC\ obtained by adding appropriate 
number of Cohen reals to a ground model, e.g.\ of \ZFC\ $+$ \CH. They 
showed that these axioms imply various assertions which were previously 
known to hold 
in Cohen models. In the present paper, we consider generalizations 
of these axioms in terms of coloring of ordinals by reals --- Homogeneity 
Principle (\HP) in section \ref{S3} and Injectivity Principle (\IP) in 
section \ref{S4}. The 
formulation of \HP\ resembles that of the Open Coloring Axiom (\OCA). 
Though it will be shown that both of our axioms are inconsistent with 
\OCA\ 
(see \Lemmaof{L2},\,\assertof{3} and \Lemmaof{L4.2},\,\assertof{3}). In 
section \ref{S5}, we 
prove that our axioms (and hence also the axioms of 
Juh\'asz, Soukup and Szentmikl\'ossy) hold not only in Cohen models but 
also in many other models of set-theory which can be quite different from 
the Cohen models from the point of view of infinitary combinatorics. 

\section{Preliminaries}\label{S2}
For a cardinal $\kappa$, $\calH(\kappa)$ denotes the set of all sets of 
hereditary cardinality $<\kappa$. In particular 
$\calH(\aleph_1)$ is the set of hereditarily countable sets. We say that a 
subset $X$ of $\calH(\aleph_1)$ is {\em definable} if there is 
a formula $\varphi(x)$ of the first order logic
in the language of 
set-theory with parameters from $\calH(\aleph_1)$ \st\ 
\[X=\setof{a\in\calH(\aleph_1)}%
	{(\calH(\aleph_1),\in)\models\varphi[a]}. \]\noindent
Note that every Borel, analytic or projective sets are definable in the 
sense as above. 

For a mapping $f$ and $X\subseteq\dom(f)$, $f``X$ denotes the image of 
$X$ by $f$. 

Following \cite{soukup-etal}, we write for a set $X$, 
$n\in\omega$ and $X_0$\tenten $X_{n-1}\subseteq X$\,:
\[ \begin{array}{c}
	(X)^n=\setof{\pairof{x_0\tenten x_{n-1}}\in X^n}{%
			x_i\not=x_j\mbox{ for all }i<j<n};\\[\jot]
	(X)^{<\aleph_0}
		=\dotbigcup_{n\in\omega}(X)^n;\\[\jot]
	(X_0\tenten X_{n-1})
		=\setof{\pairof{x_0\tenten x_{n-1}}\in (X)^n}{%
			x_0\in X_0\tenten x_{n-1}\in X_{n-1}}.
   \end{array}
\]\noindent
Note that $\psetof{\omega}$, $(\psetof{\omega})^n$, 
$(\psetof{\omega})^{<\aleph_0}$ etc.\ are definable in our sense. 
Hence we may also speak of definable subsets of these sets. 

For the consistency proof in section \ref{S5}, we use the following 
enhancement of the well-known $\Delta$-system lemma, due to Erd\H os and 
Rado. 
This theorem was already used in Juh\'asz, Soukup and Szentmikl\'ossy 
\cite{soukup-etal} and its detailed proof is also given there. 
Nevertheless, we present here another proof of the theorem 
using the method of elementary submodels. 
The proof is essentially the same as the proof of usual 
$\Delta$-system lemma given in \cite{dow}.
\begin{Thm}{\rm (Erd\H os-Rado)} \label{erdoes-rado}
Suppose that $\kappa$ and $\lambda$ are cardinals \st\ $\mu^\lambda<\kappa$ 
for all $\mu<\kappa$ and $\kappa$ is regular. If 
$\seqof{x_\alpha}{\alpha\in S}$ is a sequence of sets of cardinality 
$\leq\lambda$ for a stationary subset $S$ of 
$\setof{\alpha<\kappa}{\cof(\alpha)>\lambda}$, then there is a 
stationary $I\subseteq S$ \st\ 
$\setof{x_\alpha}{\alpha\in I}$ forms a $\Delta$-system. 
\end{Thm}
\prf
\Wolog\ we may assume that $x_\alpha\subseteq\kappa$ for every 
$\alpha\in S$. 
By Fodor's lemma and by assumption on $\kappa$ and $\lambda$, we may 
further assume that there 
is a stationary $I'\subseteq S$ \st\ $\sup(x_\alpha)\geq\alpha$ for 
every $\alpha\in I'$. 

Let $\chi$ be sufficiently large. By the 
assumptions on $\kappa$ and $\lambda$ there is $M\prec\calH(\chi)$ \st\ 
$\seqof{x_\alpha}{\alpha\in S}$, $I'$,$\kappa$, $\lambda\in M$; 
$\cardof{M}<\kappa$; $[M]^{\leq\lambda}\subseteq M$ and 
$M\cap\kappa\in I'$. Let $\alpha_1=M\cap\kappa$ and 
$r=x_{\alpha_1}\cap M$. Then $r\in M$. Hence $\alpha_0=(\sup r)+1$ is an 
element fo $M$. 

For any $\alpha_0\leq\alpha^*<\alpha_1$, we have 
$\calH(\chi)\models\alpha^*<\alpha_1\land x_{\alpha_1}\cap\alpha^*=r$. 
Hence 
$\calH(\chi)\models
	\exists\beta<\kappa(\alpha^*<\beta\land x_\beta\cap\alpha^*=r)$. 
By elementarity, it follows that 
\[M\models
	\exists\beta<\kappa(\alpha^*<\beta\land x_\beta\cap\alpha^*=r). 
\]\noindent 
Since $\alpha_1=\kappa\cap M$, it follows that 
\[M\models\forall\alpha<\kappa
	(\alpha_0\leq\alpha\rightarrow
	\exists\beta<\kappa\,(\alpha<\beta\land x_\beta\cap\alpha=r)).\]\noindent 
Again by elementarity the same sentence also holds in $\calH(\chi)$ and 
hence in $V$ as $\chi$ was taken sufficiently large. 

Now let $\xi_\alpha$, $\alpha<\kappa$ be the increasing sequence of 
ordinal numbers defined inductively by: 
$\xi_0=\mbox{the minimal }\xi\in I'\mbox{ \st\ }x_\xi\cap\alpha_0=r$
and for $\alpha>0$
\[ \begin{array}{r@{\;}l}
\xi_\alpha=&\mbox{the minimal }\xi\in I'\mbox{ \st\ }
		\xi>\sup(\bigcup\setof{x_{\xi_\beta}}{\beta<\alpha})\mbox{ and }\\
	&x_\xi\cap\sup(\bigcup\setof{x_{\xi_\beta}}{\beta<\alpha})=r.    
   \end{array}
\]\noindent
Let 
$I=\setof{\xi_\alpha}{\alpha<\kappa}$. Then $I\subseteq I'$. 
Since the definition above can be performed by using only parameters from 
$M$, we have $I\in M$. Clearly 
$\setof{x_\alpha}{\alpha\in I}$ forms a $\Delta$-system with the root 
$r$ and $I\in M$. Hence it is enough to show that $I$ is stationary. 

First note that $\alpha_1=\xi_{\alpha_1}$ and hence $\alpha_1\in I$. If 
$C\in M$ is a club subset of 
$\kappa$, then $\alpha_1\in C$. Hence $C\cap I\not=\emptyset$. By 
elementarity it follows that $M\models C\cap I\not=\emptyset$. Thus 
we have $M\models\xmbox{``$I$ is sataionary''}$. Again by elementarity 
and by the choice of $\chi$ being sufficiently large, it follows that 
$I$ is really stationary.
\qedofThm\qedskip

\memo{The condition $\mu^\lambda<\kappa$ for all $\mu<\kappa$ is a 
necessary one: suppose that $\kappa=\aleph_{\omega+1}$ and 
$\lambda=\aleph_0$. Let 
$S=\setof{\alpha<\kappa}{\cof(\alpha)=\omega}$. For $\alpha\in S$ let 
$x_\alpha$ be a cofinal subset of $\alpha$ of order type $\omega$. Then 
there is no staionary $I\subseteq S$ \st\ 
$\setof{x_\alpha}{\alpha\in I}$ build a $\Delta$-system. }
Our terminology on forcing is standard and  mainly based on K.\ 
Kunen \cite{kunen}. 
We assume that the partial orderings (or p.o.-sets, as in \cite{kunen}) 
$P$ we consider here have always the largest element $\bbbone_P$. 
For a \po\ $P$, we assume that the $P$-names are defined just as 
in \cite{kunen} but use Baumgartner's notation for $P$-names so that 
$P$-names are 
represented by dotted alphabets like $\dot{x}$, $\dot{y}$, $\dot{z}$ 
etc.\ while 
the standard names of ground model sets $x$, $y$, $z$ etc.\ are denoted by 
$\check{x}$, $\check{y}$, $\check{z}$ etc.\ or, if it is clear from the 
context, simply 
by the same symbols $x$, $y$, $z$ etc. 
Unlike \cite{kunen}, we take here the na\"\i ve (and strictly speaking 
incorrect) approach that our ground 
model is simply the set-theoretic universe $V$. 
With {\em subordering} $Q$ of $P$ we 
mean 
$Q\subseteq P$ with the ordering induced from the ordering of $P$ \st\ 
$\bbbone_P\in Q$ and hence $\bbbone_Q=\bbbone_P$; and that, for any 
$q$, $r\in Q$, if $q$ and $r$ are incompatible in $Q$ then they are 
incompatible in $P$ as well. A subordering $Q$ of a 
\po\ $P$ is called 
a {\em regular subordering} of $P$ if every maximal antichain in $Q$ is 
also a maximal antichain in $P$. If $Q$ is a regular subordering of 
$P$. As well-known, if $G$ is a $P$-generic filter over $V$ and 
$Q$ is a regular subordering of $P$, then $G\cap Q$ is $Q$-generic over $V$. 

As in \cite{kunen}, $\Fn(X,2)$ denotes 
the \po\ for adding $\cardof{X}$ many Cohen reals, i.e.\ 
$\setof{p}{\dom(p)\subseteq X,\,\range(p)\subseteq 2,\,\cardof{p}<\aleph_0}$ 
ordered by $p\leq q$ \Equivto\ $p\supseteq q$. $\dotG_P$ or simply 
$\dotG$ is the standard name of $P$ generic set which may be defined as 
$\setof{\pairof{p,\check{p}}}{p\in P}$. 

For a \po\ $Q$ and an index set $I$, the pseudo product $\Pi^*_I Q$ of 
$Q$ over $I$ was introduced in S.\ Fuchino, S.\ Shelah and L.\ Soukup 
\cite{club} as:
\[ \Pi^*_I Q=\setof{p\in\fnsp{I}{Q}}{%
		\cardof{\setof{i\in I}{p(i)\not=\bbbone_Q}}\leq\aleph_0}
\]\noindent
with the ordering
\[ p\leq q \equivto p(i)\leq q(i)\mbox{ for all }i\in I
\mbox{ and }\setof{i\in I}{p(i)\lessnoneq q(i)\lessnoneq \bbbone_Q}
\mbox{ is finite.}
\]\noindent
Recall that 
$(\clubsuit)$ is the following weakening of $\Diamond$-principle: 
\begin{assertion}{$(\clubsuit)$:}\it
There exists a sequence 
$(x_\gamma)_{\gamma\in Lim(\omega_1)}$ of countable 
subsets of $\omega_1$ \st\ for every $\gamma\in Lim(\omega_1)$, $x_\gamma$ is 
a cofinal subset of $\gamma$ with $otp(x_\gamma)=\omega$ and for every 
$y\in[\omega_1]^{\aleph_1}$ there is 
$\gamma\in Lim(\omega_1)$ \st\ $x_\gamma\subseteq y$.
\end{assertion}

Here, $Lim(\omega_1)$ denotes the set of all limit ordinals $<\omega_1$. 

In general, $(\clubsuit)$ does not hold in a Cohen model. This is because 
$\MA(Cohen)$ --- Martin's axiom restricted to the \pos\ of the form 
$\Fn(\kappa,2)$ --- can hold in such a model and, as easily seen, 
$\MA(Cohen)$ implies the 
negation of $(\clubsuit)$. However, the pseudo product of the form 
$\Pi^*_\kappa\Fn(\omega,2)$ for uncountable $\kappa$ forces $(\clubsuit)$. 
For more details and the proof of the following, see \cite{club} (for 
\assertof{2}, see also \Lemmaof{L6.5}). 

\begin{Lemma}{\rm (S.\ Fuchino, S.\ Shelah and L.\ Soukup \cite{club})} 
\label{544} 
Let $P=\Pi^*_\kappa\Fn(\omega,2)$. Then: 

\assert{1} $P$ is proper. \assertskip

\assert{2} $P$ satisfies the $(2^{\aleph_0})^+$-cc. \assertskip

\assert{3} $\forces{P}{(\clubsuit)}$. \assertskip

\assert{4} If the ground model satisfies the \CH\ and 
$\kappa^{\aleph_0}=\kappa$, then we have
\medskip\\
\mbox{}\hfill$\displaystyle\forces{P}{\MA(countable)}$.\hfill\mbox{\qed}
\end{Lemma}
\section{Homogeneity Principle}\label{S3}
One of our axioms is 
formulated as follows. First, for a cardinal $\kappa$, let 
$\HP(\kappa)$ be the following assertion:
\begin{assertion}{$\HP(\kappa)$:} \em For any 
$\mapping{f}{\kappa}{\psetof{\omega}}$ and any definable 
$A\subseteq(\psetof{\omega})^{<\aleph_0}$, either
\begin{assertion}{\phantom{\assert{h0}}}%
\mbox{}%
\lassert{h0} there is a stationary $D\subseteq\kappa$ \st\ 
$(f``D)^{<\aleph_0}\subseteq A$; or \smallskip\\
\lassert{h1} there are $n\in\omega$, $n>0$ and stationary $D_0$\tenten 
$D_{n-1}\subseteq\kappa$ \st\ 
$(f``D_0\tenten f``D_{n-1})\subseteq(\psetof{\omega})^n\setminus A$. 
\end{assertion}
\end{assertion}
By definition of definability, it is clear that $\HP(\kappa)$ remains the 
same when $\psetof{\omega}$ in this definition is replaced by one of 
$\fnsp{\omega}{\omega}$, $\fnsp{\omega}{2}$, $\reals$, $\reals^2$, $\cdots$. 

For $k\geq 1$, let $\HP_k(\kappa)$ be the assertion: 
\begin{assertion}{$\HP_k(\kappa)$:} \em For any 
$\mapping{f}{\kappa}{\psetof{\omega}}$ and any definable 
$A\subseteq(\psetof{\omega})^k$, either
\begin{assertion}{\phantom{\assert{h0'}}}%
\mbox{}%
\lassert{h0'} there is a stationary $D\subseteq\kappa$ \st\ 
$(f``D)^k\subseteq A$; or \smallskip\\
\lassert{h1'} there are stationary $D_0$\tenten 
$D_{k-1}\subseteq\kappa$ \st\\
$(f``D_0\tenten f``D_{k-1})\subseteq(\psetof{\omega})^k\setminus A$. 
\end{assertion}
\end{assertion}
\begin{Lemma}\label{L1}

\assert{0} $\HP_1(\kappa)$ holds for every regular $\kappa$.\assertskip

\assert{1} If $1\leq k<k'<\omega$, then $\HP(\kappa)$ implies 
$\HP_{k'}(\kappa)$ and $\HP_{k'}(\kappa)$ implies 
$\HP_k(\kappa)$.\assertskip

\assert{2} If $\cof(\kappa)=\lambda$ and $\HP(\lambda)$, then 
$\HP(\kappa)$ also holds. \assertskip

\assert{3} $\HP(\kappa)$ holds for any regular $\kappa>2^{\aleph_0}$. 
\assertskip

\assert{4} If\/ $\HP_2(\kappa)$, then there is no subset $X$ of 
$\psetof{\omega}$ with a definable $A\subseteq\psetof{\omega}^2$ \st\ 
$A\cap X^2$ is a well-ordering on $X$ of order-type $\kappa$. In 
particular, $\kappa$ is not embeddable into 
$(\psetof{\omega}/fin,\subseteq^*)$. 
\assertskip

\assert{5} $\neg\HP_2(\aleph_1)$. Hence by \assertof{1}, we have 
$\neg\HP_k(\aleph_1)$ for every $k\geq 2$ and $\neg\HP(\aleph_1)$. 
\end{Lemma}
\prf\assertof{0}\,: Trivial. 

\assertof{1}\,: To show that $\HP(\kappa)$ implies $\HP_{k'}(\kappa)$, let 
$\mapping{f}{\kappa}{\psetof{\omega}}$ and 
$A'\subseteq(\psetof{\omega})^{k}$ be definable. Let 
\[ \begin{array}{r@{\;}l}
A=&(\psetof{\omega})^1\cup\cdots\cup(\psetof{\omega})^{k'-1}\cup\\[\jot]
	&{\bigcup}_{l\geq k'}
		\setof{\pairof{x_0\tenten x_{l-1}}\in(\psetof{\omega})^l}{%
			\pairof{x_0\tenten x_{k-1}}\in A'}.   
   \end{array}
\]\noindent
Then if \assertof{h0} holds 
for $f$ and $A$ above, we can easily see that \assertof{h0'} of 
$HP_{k'}(\kappa)$ holds for $f$ and $A'$. Similarly for \assertof{h1}. 
The implication 
of $\HP_k(\kappa)$ form $\HP_{k'}(\kappa)$ can be shown similarly. 

\assertof{2}\,: Easy.

\assertof{3}\,: 
Let $\kappa>2^{\aleph_0}$ be a regular cardinal. Suppose that 
$\mapping{f}{\kappa}{\psetof{\omega}}$ and $A$ are as in the definition of 
$\HP(\kappa)$. 
Then there is a stationary 
$D\subseteq\kappa$ \st\ $f\restr D$ is constant. If \assertof{h0} in the 
definition of $\HP(\kappa)$ does not hold then we must have 
$(f``D)^{1}\cap A=\emptyset$ since $(f``D)^{n}=\emptyset$ for $n>1$. 
Hence \assertof{h1} holds with $n=1$ and 
$D_0=D$. 

\assertof{4}\,: Suppose that $X\subseteq\psetof{\omega}$, 
$A\subseteq(\psetof{\omega})^2$ is definable and $A\cap (X)^2$ is a 
well-ordering on $X$ of order type $\kappa$. 
If $\mapping{f}{\kappa}{\psetof{\omega}}$ is the mapping sending $\alpha$ to 
the $\alpha$'th element of $X$ \wrt\ $A$ then 
none of conditions \assertof{h0'} and \assertof{h1'} in the definition of 
$\HP_2(\kappa)$ can hold: for any stationary $D\subseteq\kappa$, if 
$\alpha$, $\beta\in D$ are \st\ $\beta<\alpha$, then we have 
$\pairof{\alpha,\beta}\not\in A$; on the other hand, if $D_0$, 
$D_1$ are stationary subsets of $\kappa$, then, taking 
$\pairof{\alpha,\beta}\in (D_0,D_1)$ \st\ $\alpha<\beta$, we have 
$\pairof{\alpha,\beta}\in A$. 

\assertof{5} follows immediately from \assertof{4} since there is a subset of 
$\psetof{\omega}/fin$ of order-type $\omega_1$ \wrt\ $\subseteq^*$. 
\qedofLemma\qedskip

The following is still open:
\begin{Problem}
Does $\HP_2(\kappa)$ imply $\HP_k(\kappa)$ for $k\geq 3$? Does it imply 
$\HP(\kappa)$? 
\end{Problem}

Now we define the {\em Homogeneity Principle} (\HP) as the assertion:
\begin{assertion}{$\HP$:}
$\HP(\kappa)$ holds for every regular $\kappa>\aleph_1$. 
\end{assertion}
\begin{Lemma}\label{L2}

\assert{1} \CH\ implies \HP.\smallskip

\assert{2} $\HP_2(\aleph_2)$ and hence $\HP$ implies 
${\bf b}=\aleph_1$. \smallskip 

\assert{3} \OCA\ implies $\neg\HP_2(\aleph_2)$ and hence $\neg\HP$. \smallskip

\assert{4} $\HP_2(\aleph_2)$ and hence $\HP$ implies that 
$\psetof{\omega_1}/countable$ is not embeddable into $\psetof{\omega}/fin$. 
\end{Lemma}
\prf \assertof{1} follows from \Lemmaof{L1},\,\assertof{3}. 

\assertof{2}\,: By \Lemmaof{L1},\,\assertof{4}, it follows  from 
$\HP_2(\aleph_2)$ that ${\bf b}<\aleph_2$. 

\assertof{3}\,: \OCA\ implies ${\bf b}=\aleph_2$ (see e.g.\ 
\cite{bekkali}). By \assertof{2}, this implies $\neg\HP_2(\aleph_2)$. 

\assertof{4}\,: It is known that $\omega_2$ is embeddable into 
$\psetof{\omega_1}/countable$ (see e.g.\ \cite{FuKoSh}, p48). Hence the 
assertion follows from \Lemmaof{L1},\,\assertof{4}. 

\qedofLemma\qedskip

Let $\HP^+(\kappa)$ be the assertion which is obtained when the 
assumption on $A$ being definable is dropped from
the definition of $\HP(\kappa)$ and $\HP^+$ be the corresponding axiom. 
The proof of \Lemmaof{L1},\,\assertof{4} shows that 
$\HP^+(\kappa)$ cannot hold 
for any $\kappa\leq2^{\aleph_0}$ under \ZFC. 
On the other hand, the proof of 
\Lemmaof{L1},\,\assertof{3} works also for this stronger form of the 
homogeneity principle. Hence we see that $\HP^+$ is just equivalent to 
$\CH$. 

One of the principles considered by Juh\'asz, Soukup and Szentmikl\'ossy 
\cite{soukup-etal} is the following $C^s(\kappa)$ for a cardinal $\kappa$. 

\begin{assertion}{$C^s(\kappa)$:}\it
For any $T\subseteq\omega^{<\omega}$ and any matrix 
$\seqof{a(\alpha,n)}{\alpha<\kappa,n\in\omega}$ of subsets of 
$\omega$, either 
\begin{assertion}{\phantom{\assert{c0}}}%
\mbox{}\lassert{c0} there is a stationary $S\subseteq\kappa$ \st\ 
for each $t\in T$ and $s\in(S)^{|t|}$ we have 
$\bigcap_{i<|t|}a(s(i),t(i))\not=\emptyset$;
or\smallskip\\
\lassert{c1} there are $t\in T$ and stationary $D_0$\tenten 
$D_{|t|-1}\subseteq\kappa$ \st\ for every 
$s\in (D_0\tenten D_{|t|-1})$ we have 
$\bigcap_{i<|t|}a(s(i),t(i))=\emptyset$. 
\end{assertion}
\end{assertion}
A sort of the dual of $C^s(\kappa)$ is called $\hat{C}^s(\kappa)$ by 
Juh\'asz, Soukup and Szentmikl\'ossy \cite{soukup-etal}:
\begin{assertion}{$\hat{C}^s(\kappa)$:}\it
For any $T\subseteq\omega^{<\omega}$ and any matrix 
$\seqof{a(\alpha,n)}{\alpha<\kappa,n\in\omega}$ of subsets of 
$\omega$, either 
\begin{assertion}{\phantom{\assert{c0}}}%
\mbox{}\lassert{\^c0} there is a stationary $S\subseteq\kappa$ \st\ 
for each $t\in T$ and $s\in(S)^{|t|}$ we have 
$\cardof{\bigcap_{i<|t|}a(s(i),t(i))}<\aleph_0$;
or\smallskip\\
\lassert{\^c1} there are $t\in T$ and stationary $D_0$\tenten 
$D_{|t|-1}\subseteq\kappa$ \st\ for every 
$s\in (D_0\tenten D_{|t|-1})$ we have 
$\cardof{\bigcap_{i<|t|}a(s(i),t(i))}=\aleph_0$. 
\end{assertion}
\end{assertion}
The weakenings of these assertions $C(\kappa)$ and $\hat{C}(\kappa)$ are also 
considered in \cite{soukup-etal}, which are obtained when we replace 
``stationary \ldots'' by ``\ldots\ of cardinality $\kappa$'' in the 
definition of 
$C^s(\kappa)$ and $\hat{C}^s(\kappa)$ respectively. 
Note that we could also consider the version of $\HP(\kappa)$ and 
$\HP$ obtained by replacing stationarity by cardinality $\kappa$; it is 
easily checked that Lemmas \ref{L1}, \ref{L2} still hold under this weaker 
version of Homogeneity Principle.

\begin{Lemma}\label{L3}
For every $\kappa$, $\HP(\kappa)$ implies $C^s(\kappa)$ and 
$\hat{C}^s(\kappa)$. 
\end{Lemma}
\prf
We show that $\HP(\kappa)$ implies $C^s(\kappa)$. The assertion for 
$\hat{C}^s(\kappa)$ can be proved similarly. 

First, with the same argument as in the proof of \Lemmaof{L1},\assertof{3}, 
we can prove easily that $C^s(\kappa)$ (and also $\hat{C}^s(\kappa)$) 
holds for every regular $\kappa>2^{\aleph_0}$. Hence \wolog\ we may 
assume that $\kappa\leq2^{\aleph_0}$. 

Let 
$\seqof{t_i}{i\in\omega}$ be an enumeration of $\omega^{<\omega}$ \st\ 
$\cardof{t_i}\leq i$ for every $i\in\omega$. 
Let $\mapping{\iota}{\psetof{\omega}}{\psetof{\omega}^{\aleph_0}}$ be some 
definable bijection. For $x\in\psetof{\omega}$ and let 
$(x)_i$ denote the $i$'th component of $\iota(x)$ for each $i\in\omega$. 

Now suppose that $T\subseteq\omega^{<\omega}$ and 
$\seqof{a(\alpha,n)}{\alpha<\kappa,\,n\in\omega}$ is a matrix of subsets 
of $\omega$. Let $\mapping{g}{\kappa}{\psetof{\omega}}$ be a fixed 
inective mapping (which exists by the assumption 
$\kappa\leq 2^{\aleph_0}$) and let 
$\mapping{f}{\kappa}{\psetof{\omega}}$ be defined by 
$f(\alpha)=\iota^{-1}(\seqof{a'(\alpha,n)}{n\in\omega})$ where 
$a'(\alpha,n)$'s are defined by $a'(\alpha,0)=g(\alpha)$ and 
$a'(\alpha, n+1)=a(\alpha,n)$. Note that $f$ is 1-1 mapping because of 
the clause on $a'(\alpha,0)$. 
For each 
$i<\omega$, let 
\[ A_i=\left\{\,
\begin{array}{@{}ll}
\setof{\pairof{x_0\tenten x_{i-1}}\in(\psetof{\omega})^i}%
	{{\bigcap}_{n<\cardof{t_i}}(x_n)_{t_i(n)+1}\not=\emptyset},\quad
		&\mbox{if }t_i\in T;\\[\jot]
(\psetof{\omega})^i,
		&\mbox{otherwise}
\end{array}\right.
\]\noindent
and let $A=\dotbigcup_{i\in\omega}A_i$. Since $T$ is hereditary countable 
and $\iota$ is definable, $A$ is also 
definable. Hence, for these $f$ and $A$,  either \assertof{h0} or 
\assertof{h1} in 
the definition of $\HP(\kappa)$ holds. 

If \assertof{h0}, then there is a stationary $D\subseteq\kappa$ \st\ 
$(f``D)^{<\aleph_0}\subseteq A$. We claim that this $D$ witnesses
\assertof{c0} in the definition of $C^s(\kappa)$ for these $T$ and 
$\seqof{a(\alpha,n)}{\alpha<\kappa,\,n\in\omega}$: for any $t\in T$ and 
$s\in(D)^{\cardof{t}}$, let $i\in\omega$ be \st\ $t=t_i$. Let 
$s'\in(D)^i$ be any end-extension of $s$. 
Since $f$ is 1-1, 
we have 
$\pairof{f(s'(0))\tenten f(s'(i-1))}\in (f``D)^i$. 
Hence $\pairof{f(s'(0))\tenten f(s'(i-1))}\in A_i$. By 
definition of $A_i$ this means 
$\bigcap_{n<\cardof{t}}a(s(n),t(n))\not=\emptyset$. 

If 
\assertof{h1}, then there is an $i<\omega$ and stationary 
$D_0$\tenten $D_{i-1}\subseteq\kappa$ \st\ 
$(f``D_0\tenten f``D_{i-1})\subseteq(\psetof{\omega})^i\setminus A_i$. 
Let $t=t_i$. Then $t\in T$ by definition of $A$ and
$\bigcap_{j<|t|}a(s(j),t(j))=\emptyset$ for every 
$s\in (D_0\tenten D_{|t|-1})$. Thus \assertof{c1} holds for these 
$D_0$\tenten $D_{|t|-1}$. \qedofLemma\qedskip

Using the trick of the proof above we can also easily obtain the 
following lemma which generalizes Theorem 4.1 in \cite{soukup-etal}:
\begin{Lemma}\label{L4}
Suppose that $\HP(\kappa)$. 
For any sequence $\seqof{x_\alpha}{\alpha<\kappa}$ of subsets of 
$\omega$, either there is a stationary $D\subseteq\kappa$ \st\ 
$\setof{x_\alpha}{\alpha\in D}$ is independent (i.e.\ equivalence 
classes of $x_\alpha$, $\alpha\in D$ are independent in 
$\psetof{\omega}/fin$ as a subset of the Boolean algebra); or there are 
$n\in\omega$, $\mapping{i}{n}{2}$ and stationary 
$D_0$\tenten $D_{n-1}\subseteq\kappa$ \st\ 
$\bigcap_{k<n}(x_{\alpha_k})^{i(k)}$ is finite for every 
$\pairof{\alpha_0\tenten\alpha_{n-1}}\in (D_0\tenten D_{n-1})$,  
where we write
\[ (x)^i=\left\{
\begin{array}{@{\,\,}ll}
x\quad&\mbox{if }i=1,\\[\jot]
\omega\setminus x\quad&\mbox{otherwise,}
\end{array}\right.
\]\noindent
for $x\in\psetof{\omega}$ and $i\in 2$. 
\qed
\end{Lemma}

By \Lemmaof{L3}, every assertions provable under $C^s(\kappa)$ or 
$\hat{C}^s(\kappa)$ e.g.\ such as given in \cite{soukup-etal} are also 
consequences of $\HP(\kappa)$. One of such assertions is the 
non-existence of 
$\kappa$-Lusin gap for which we shall give here a direct proof from 
$\HP(\kappa)$. 

An almost disjoint family $F\subseteq[\omega]^{\aleph_0}$ is called a 
{\em$\kappa$-Lusin gap} if, $\cardof{F}=\kappa$ and, for any 
$x\in[\omega]^{\aleph_0}$, at least one of 
$\cardof{\setof{y\in F}{y\subseteq^*x}}<\kappa$ or 
$\cardof{\setof{y\in F}{y\subseteq^*\omega\setminus x}}<\kappa$ holds. 
\begin{Lemma}{\rm (Juh\'asz, Soukup and Szentmikl\'ossy 
\cite{soukup-etal} under $C(\kappa)$)} Suppose that $\kappa$ is a regular 
cardinal and $\HP(\kappa)$ holds. 
Then there is no $\kappa$-Lusin gap. 
\end{Lemma}
\prf
Suppose that  $F\subseteq[\omega]^{\aleph_0}$ is an almost disjoint 
family of cardinality $\kappa$. 
We show that $F$ is no $\kappa$-Lusin gap. 
Let $\mapping{f}{\kappa}{F}$ be an injective mapping and let
\[ A=\setof{\pairof{x_0\tenten x_{n-1}}\in(\psetof{\omega})^{<\omega}}%
		{n=1\mbox{ or }(x_0\cap x_1)\setminus n\not=\emptyset}. 
\]\noindent
Then clearly \assertof{h0} in the definition of $\HP(\kappa)$ is 
impossible for these $f$, $A$. Hence, by $\HP(\kappa)$, there is 
$k\in\omega$, $k\geq 1$, and stationary $D_0$\tenten 
$D_{k-1}\subseteq\kappa$ \st\ 
$(f``D_0,f``D_1$\tenten $f``D_{k-1})\subseteq(\psetof{\omega})^k\setminus A$. 
By definition of $A$, this means that 
$(\bigcup_{\alpha\in D_0}f(\alpha))\cap(\bigcup_{\alpha\in D_1}f(\alpha))
	\subseteq k$. 

Let $x=\bigcup_{\alpha\in D_0}f(\alpha)$. Then 
$\setof{f(\alpha)}{\alpha\in D_0}\subseteq\setof{y\in F}{y\subseteq^* x}$ and 
$\setof{f(\alpha)}{\alpha\in D_1}\subseteq
	\setof{y\in F}{y\subseteq^* (\omega\setminus x)}$. Since $f``D_0$ and 
$f``D_1$ are both of cardinality $\kappa$, this shows that $F$ is not a 
$\kappa$-Suslin gap. 
\qedofLemma

\section{Injectivity Principle}\label{S4}
In this section we consider two axioms which are connected with the 
principles 
$F^s(\kappa)$ introduced in Juh\'asz, Soukup and Szentmikl\'ossy 
\cite{soukup-etal}. For a regular cardinal $\kappa$ and 
$\lambda\leq\kappa$, stipulate 
\begin{assertion}{$\IP(\kappa,\lambda)$:}
\em For any $\mapping{f}{\kappa}{\psetof{\omega}}$ and definable 
$\mapping{g}{(\psetof{\omega})^{<\omega}}{\psetof{\omega}}$, either
\begin{assertion}{\phantom{000}}%
\mbox{}\lassert{i0} there is a stationary set $D\subseteq\kappa$ \st\ 
$\cardof{g``(f``D)^n}<\lambda$ for every $n\in\omega$; or\smallskip\\
\mbox{}\lassert{i1} there are $n\in\omega$ and stationary sets 
$D_0$\tenten $D_{n-1}\subseteq\kappa$ \st\ for any 
$\pairof{x_0\tenten x_{n-1}}$, 
$\pairof{y_0\tenten y_{n-1}}\in (f``D_0\tenten f``D_{n-1})$, if $x_i\not= y_i$ for 
all $i<n$, then we have $g(x_0\tenten x_{n-1})\not=g(y_0\tenten y_{n-1})$. 
\end{assertion}
\end{assertion}
As for $\HP(\kappa)$, we may replace $\psetof{\omega}$ in the definition of 
$\IP(\kappa,\lambda)$ by one of 
$\fnsp{\omega}{\omega}$, $\fnsp{\omega}{2}$, $\reals$, $\reals^2$, $\cdots$. 
We can consider the weakenings $\IP_k(\kappa,\lambda)$, 
$k\in\omega$ which corresponds to $\HP_k(\kappa)$, $k\in\omega$ and 
reformulate the following lemmas according to these weakenings of 
$\IP(\kappa,\lambda)$. Also, for the most of the following, the version 
of $\IP(\kappa,\lambda)$ is enough which is obtained by replacing 
stationarity in the definition by cardinality $\kappa$. 

\begin{Lemma}\label{L4.1}
Suppose that $\kappa$ is a regular uncountable cardinal and $\lambda<\kappa$. 
\smallskip

\assert{0} If $\lambda'<\lambda$ and $\IP(\kappa,\lambda')$ then 
$\IP(\kappa,\lambda)$. \smallskip

\assert{1} If $2^{\aleph_0}<\kappa$ then $\IP(\kappa,2)$.\smallskip

\assert{2} $\IP(\kappa,\aleph_0)$ implies $2^{\aleph_0}<\kappa$.\smallskip

\assert{3} $\IP(\aleph_2,\aleph_0)$ is equivalent to \CH. \smallskip

\assert{4} $\IP(\kappa,\kappa)$ implies that there is no strictly 
increasing chain in 
$\psetof{\omega}$ of length $\kappa$ \wrt\ $\subseteq^*$. 
\smallskip

\assert{5} $\neg\IP(\aleph_1,\aleph_1)$. \smallskip
\end{Lemma}
\prf \assertof{0} is clear by definition. 

\assertof{1}\,: Let $f$ and $g$ be as in the definition of $\IP(\kappa,2)$. 
Let $D\subseteq\kappa$ be a stationary set \st\ $f$ is constant on 
$D$. Then clearly \assertof{i0} holds for this $D$. 

\assertof{2}\,: Suppose that $2^{\aleph_0}\geq\kappa$. We show that 
$\IP(\kappa,\aleph_0)$ does not hold. Let 
$\mapping{f}{\kappa}{\psetof{\omega}}$ be injective and 
$\mapping{g}{(\psetof{\omega})^{<\omega}}{\psetof{\omega}}$ be defined by 
$g(x_0)=\emptyset$ and 
$g(x_0\tenten x_{n-1})=\min\setof{n\in\omega}%
	{n\in x_0\not\leftrightarrow n\in x_1}$ for $n\geq2$.
Let $D$ be any stationary subset of $\kappa$. Then 
$\cardof{g``(f``D)^2}\geq\aleph_0$: otherwise, let $k\in\omega$ be \st\ 
$g``(f``D)^2\subseteq k$. Since $\psetof{k+1}$ is finite, there are 
$\alpha$, $\beta\in D$, $\alpha\not=\beta$ \st\ 
$f(\alpha)\restr (k+1)=f(\beta)\restr (k+1)$. But then, by definition of 
$g$, it follows that $g(f(\alpha), f(\beta))>k$. This is a contradiction. 
Thus \assertof{i0} does not hold for these $f$ and $g$. 
On the other hand, 
for arbitrary stationary subsets $D_0$\tenten $D_{n-1}$ of $\kappa$, as 
there are only countably many values of $g$, we can easily find 
$\pairof{x_0\tenten x_{n-1}}$, 
$\pairof{y_0\tenten y_{n-1}}\in (f``D_0\tenten f``D_{n-1})$  \st\ $x_i\not= y_i$ for 
all $i<n$ and $g(x_0\tenten x_{n-1})=g(y_0\tenten y_{n-1})$. Thus 
neither \assertof{i1} can hold. 

\assertof{3} follows from \assertof{1} and \assertof{2}. 

\assertof{4}\,: Suppose that $\seqof{x_\alpha}{\alpha<\kappa}$ is a 
strictly increasing chain in $\psetof{\omega}$ \wrt\ $\subseteq^*$. 
Let 
$\mapping{f}{\kappa}{\psetof{\omega}}$; $\alpha\mapsto x_\alpha$ and 
$\mapping{g}{(\psetof{\omega})^{<\omega}}{\psetof{\omega}}$ be defined by 
$f(x_0)=\emptyset$ and $f(x_0\tenten x_{n-1})=x_0\setminus x_1$ for 
$n\geq 2$. We claim that these $f$ and $g$ are counter examples for 
$\IP(\kappa,\aleph_1)$. For arbitrary stationary subset $D$ of 
$\kappa$, if $\alpha_0$, $\alpha_1$, $\beta_0$, $\beta_1\in D$ are \st\ 
$\beta_0<\alpha_0\leq\beta_1<\alpha_1$, then 
$f(\alpha_0)\setminus f(\beta_0)$ and 
$f(\alpha_1)\setminus f(\beta_1)$ are two almost disjoint infinite sets. 
It follows that $\cardof{g``(f``D)^2}=\kappa$. Hence \assertof{i0} does 
not hold for these $f$, $g$. On the other hand, if $D_0$\tenten 
$D_{n-1}$ are stationary, then 
$\setof{g(f(\alpha_0)\tenten f(\alpha_{n-1}))}%
	{\pairof{\alpha_0\tenten\alpha_{n-1}}\in (D_0\tenten D_{n-1}),\,
	\alpha_0<\alpha_1}$ 
is a subset of $[\omega]^{<\aleph_0}$ and hence countable. So we can find 
$\pairof{x_0\tenten x_{n-1}}$, 
$\pairof{y_0\tenten y_{n-1}}\in (f``D_0\tenten f``D_{n-1})$  with $x_i\not= y_i$ for 
all $i<n$ \st\ $g(x_0\tenten x_{n-1})=g(y_0\tenten y_{n-1})$. Thus 
\assertof{i1} cannot hold for these $f$, $g$. 

\assertof{5}\,: There exists a strictly increasing chain in 
$\psetof{\omega}$ \wrt\ $\subseteq^*$ of length $\omega_1$. Hence, by 
\assertof{4}, we have $\neg\IP(\aleph_1,\aleph_1)$. 
\qedofLemma\qedskip

Now, similarly to the Homogeneity Principle, we define the Injectivity 
Principle (\IP) as the assertion:
\begin{assertion}{\IP:}
\em
$\IP(\kappa,\kappa)$ holds for every regular cardinal $\kappa\geq\aleph_2$. 
\end{assertion}
A natural strengthening of \IP\ would be the following:
\begin{assertion}{$\IP^+$:}
\em
$\IP(\kappa,\aleph_1)$ holds for every regular cardinal $\kappa\geq\aleph_2$. 
\end{assertion}
The following lemma can be now proved similarly to \Lemmaof{L2}: 

\begin{Lemma}\label{L4.2}

\assert{1} \CH\ implies \IP.\smallskip

\assert{2} $\IP(\aleph_2,\aleph_2)$ and hence $\IP$ implies 
${\bf b}=\aleph_1$. \smallskip 

\assert{3} \OCA\ implies $\neg\IP(\aleph_2,\aleph_2)$ and hence 
$\neg\IP$. \smallskip

\assert{4} $\IP(\aleph_2,\aleph_2)$ and hence $\IP$ implies 
$\psetof{\omega_1}/countable$ is not embeddable into $\psetof{\omega}/fin$. 
\qed
\end{Lemma}

For a regular cardinal $\kappa$, the principle $F^s(\kappa)$ is defined in 
\cite{soukup-etal} by:
\begin{assertion}{$F^s(\kappa)$:}\it
For any $T\subseteq\omega^{<\omega}$ and any matrix 
$\seqof{a(\alpha,n)}{\alpha<\kappa,n\in\omega}$ of subsets of 
$\omega$, either 
\begin{assertion}{\phantom{\assert{f0}}}%
\mbox{}\lassert{f0} there is a stationary $S\subseteq\kappa$ \st\
\[ \cardof{\setof{{\bigcap}_{i<|t|}a(s(i),t(i))}%
		{t\in T\mbox{ and }s\in (S)^{|t|}}}\leq\aleph_0; or
\]\noindent
\lassert{f1} there are $t\in T$ and stationary $D_0$\tenten 
$D_{|t|-1}\subseteq\kappa$ \st\ for every 
$s_0, s_1\in (D_0\tenten D_{|t|-1})$, if $s_0(i)\not= s_1(i)$ for all 
$i<|t|$, then 
$\bigcap_{i<|t|}a(s_0(i),t(i))\not=\bigcap_{i<|t|}a(s_1(i),t(i))$. 
\end{assertion}
\end{assertion}
Using the trick of the proof of \Lemmaof{L3}, the following is readily seen: 
\begin{Lemma}
$\IP(\kappa,\aleph_1)$ implies $F^s(\kappa)$. Thus $\IP^+$ implies 
$F^s(\kappa)$ for every regular $\kappa\geq\aleph_2$. 
\qed
\end{Lemma}
\section{Nice side-by-side products}\label{S5}
In Juh\'asz, Soukup and Szentmikl\'ossy \cite{soukup-etal}, it is proved 
that, if $\kappa$ is an $\aleph_0$-inaccessible cardinal 
in the ground model (for definition see below), then 
$C^s(\kappa)$, $\hat{C}^s(\kappa)$ and $F^s(\kappa)$ still hold in the model 
obtained by adding arbitrary number of Cohen reals. Practically the same 
proof also shows that $\HP(\kappa)$ and 
$\IP(\kappa)$ hold in such a model. Analyzing this proof, it is 
easily seen that the idea of the proof works not only for the \pos\ of 
the form $\Fn(\lambda,2)$ in the case of Cohen models but also for a much 
broader class of \pos. 

To formulate this fact in a general form, we introduce the following 
class of \pos: for a cardinal $\kappa$ and a \po\ $Q$, we 
call a \po\ $P$ a {\em nice side-by-side $\kappa$-product of $Q$} if 
there are suborderings $Q_\alpha$, $\alpha<\kappa$ of $P$, 
order-isomorphisms $\mapping{\iota_\alpha}{Q}{Q_\alpha}$ and order 
preserving ``projections'' $\mapping{\pi_\alpha}{P}{P_\alpha}$ for all 
$\alpha<\kappa$ 
with the following properties \assertof{n1} -- 
\assertof{n4}:\medskip 

\assert{n1} For every $p\in P$ 
$\setof{\alpha<\kappa}{\pi_\alpha(p)\not=\bbbone_P}$ is countable and 
$p=\inf\setof{\pi_\alpha(p)}{\alpha<\kappa}$. 
We denote with $\supp(p)$ the countable set 
$\setof{\alpha<\kappa}{\pi_\alpha(p)\not=\bbbone_P}$. 
\smallskip

\assert{n2} For any $p\in P$ and $X\subseteq\kappa$, 
$\inf\setof{\pi_\alpha(p)}{\alpha\in X}$ exists and 
$\supp(p\restr X)\subseteq X$. We denote this element of 
$P$ with $p\restr X$. 
\smallskip

\assert{n3} For $X\subseteq\kappa$, let 
$P_X=\setof{p\in P}{\supp(p)\subseteq X}$. Then 
the mapping 
$P_X\times P_{\kappa\setminus X}\ni \pairof{p,q}
	\mapsto\inf\ssetof{p,q}\in P$ 
is well-defined and is an order-isomorphism.\smallskip

\assert{n4} Every bijection $\mapping{j}{\kappa}{\kappa}$ induces an 
automorphism $\mapping{\tilde{j}}{P}{P}$ \st\ 
$\tilde{j}\restr Q_\alpha=\iota_{j(\alpha)}\circ\iota_\alpha^{-1}$. 
\medskip\\ 

Note that the usual 
finite and countable support side-by-side product of $Q$ as well as 
$\Pi^*_\kappa Q$ satisfy the conditions above. Also, note that the Cohen 
forcing $\Fn(\kappa,2)$ is isomorphic to the finite support 
$\kappa$-product of $\Fn(\omega,2)$. 

\begin{Lemma}\label{L4.5} Suppose that $P$ is a nice side-by-side 
$\kappa$ product of a \po\ $Q$ and $X\subseteq P$. Then, for any $p$, 
$q\in P$: \assertskip

\assert{0} $P_X$ is a regular subordering of $P$.\assertskip

\assert{1} 
The mapping 
$P\ni p\mapsto
	\pairof{p\restr X,p\restr(\kappa\setminus X)}
		\in P_X\times P_{\kappa\setminus X}$ 
is the reverse mapping of the isomorphism given in 
\assertof{n3}. In particular, if $p\leq q$, then 
$p\restr X\leq q\restr X$. 
\assertskip

\assert{2} If $X=\supp(p)\cap\supp(q)$ and $p\restr X\leq q\restr X$, then 
$r=\inf\setof{p(\alpha), q(\alpha)}{\alpha\in\kappa}$ exists. We shall write 
$p\land q$ to denote such $r$.\assertskip

\assert{3} Let $\lambda=\cardof{X}$. Then  $P_X$ is a nice 
side-by-side $\lambda$-product of $Q$. 

\end{Lemma}
\prf \assertof{0}, \assertof{1} and \assertof{2} follows immediately from 
\assertof{n2} and \assertof{n3}. 
For \assertof{3}, enumerate $Q_\alpha$, $\alpha\in X$ by $\lambda$, 
say as $\seqof{Q'_\alpha}{\alpha<\lambda}$. Then it 
is easy to see that this sequence witnesses that $P\restr X$ is a nice 
side-by-side $\lambda$-product of $Q$. \qedofLemma\qedskip

In \cite{soukup-etal}, a cardinal number $\kappa$ is called 
$\lambda$-inaccessible if $\kappa$ is regular and $\mu^\lambda<\kappa$ 
for every $\mu<\kappa$. Under $\GCH$, every regular 
cardinal $\kappa\geq\aleph_2$ is $\aleph_0$-inaccessible provided that 
$\kappa$ is not a successor of a cardinal of cofinality $\omega$. 

\begin{Thm}\label{main-thm-1}
Let $\kappa$ be an $\aleph_0$-inaccessible cardinal. 
Suppose that $P$ is a proper nice side-by-side $\kappa$-product of a \po\ 
$Q$ \st\ $\cardof{Q}<\kappa$. Then $\forces{P}{\HP(\kappa)}$.
\end{Thm}

Toward the proof of \Thmabove, we prove first the following lemmas. 

For a \po\ $P$, a $P$-name $\dotx$ of a subset of $\omega$ is called a 
{\it nice $P$-name} if there is an antichain $A^\dotx_n$ in $P$ for each 
$n\in\omega$ \st\ 
$\dotx=\bigcup_{n\in\omega}\setof{\pairof{p,\check n}}{p\in A^\dotx_n}$. 
It is easy to see that, for each $P$-name $\dotx$ of a subset of 
$\omega$, there is a nice 
$P$-name $\dotx'$ \st\ $\forces{P}{\dotx=\dotx'}$. We say that a nice 
$P$-name 
$\dotx$ is {\em slim} if each $A^\dotx_n$ as above is countable. 
If $\dotx$ is a nice $P$-name of a subset of $\omega$ with $A^\dotx_n$, 
$n\in\omega$ as above and $P$ is a nice 
side-by-side $\kappa$-product of $Q$ for some $\kappa$ and $Q$, then we 
define 
$\supp(\dotx)=
	\bigcup\setof{\supp(p)}{p\in A^\dotx_n\xmbox{ for some }n\in\omega}$. 
Note that, if $\dotx$ is slim, then $\supp(\dotx)$ is countable. 
An automorphism $\tilde j$ of $P$ induces canonically a bijective class 
function on $P$-names which is denoted also by $\tilde j$. For a nice 
$P$-name $\dotx$, this function maps $\dotx$ to 
$\tilde{j}(\dotx)=
	\setof{\pairof{\tilde{j}(p),\check n}}{\pairof{p,\check n}\in \dotx}$. 

\begin{Lemma}\label{L5}
If $P$ is proper  and $p\in P$ then, 
for any name $\dotx$ of a subset of 
$\omega$, there are $q\leq p$ and a slim $P$-name $\dotx'$ \st\ 
$q\forces{P}{\dotx=\dotx'}$. 
\end{Lemma}
\prf By the remark above we may assume \wolog\ that $\dotx$ is a nice 
$P$-name. Let $A^\dotx_n$, $n\in\omega$ be as above and 
$\doty$ be a $P$-name \st\ 
$\forces{P}{\doty=\setof{s\in P}%
	{s\in(\bigcup_{n\in\omega}A^\dotx_n)\cap\dotG}}$. Then we have 
$\forces{P}{\doty\xmbox{ is a countable subset of }P}$. As $P$ is proper 
there exist $q\leq p$ and countable $y\subseteq P$ \st\ 
$q\forces{P}{\doty\subseteq y}$. Let 
$\dotx'=\bigcup_{n\in\omega}
	\setof{\pairof{s,\check n}}{s\in A^\dotx_n\cap y}$. It is 
then easy to see that these 
$q$ and $\dotx'$ are as desired. \qedofLemma

\begin{Lemma}\label{L6} Suppose that $\kappa$ is 
$\aleph_0$-inaccessible and $P$ is a nice side-by-side 
$\kappa$-product of a \po\ $Q$ with $\cardof{Q}<\kappa$. 
If $S\subseteq\setof{\alpha<\kappa}{\cof(\alpha)>\omega}$ is stationary and 
$\seqof{\dotx_\alpha}{\alpha\in S}$ is a sequence of $P$-names of subsets 
of $\omega$, then for any $p\in P$ there are a stationary 
$S^*\subseteq S$, a sequence 
$\seqof{\dotx'_\alpha}{\alpha\in S^*}$ of slim $P$-names and a 
sequence $\seqof{p_\alpha}{\alpha\in S^*}$ of elements of $P$ \st\medskip 

\assert{a} $p_\alpha\leq p$ and 
$p_\alpha\forces{P}{\dotx_\alpha=\dotx'_\alpha}$ for every $\alpha\in S^*$;

\assert{b} $d_\alpha=\supp(p_\alpha)\cup\supp(\dotx'_\alpha)$, 
$\alpha\in S^*$ are all of the same cardinality and form 
a $\Delta$-system with a root $r$;

\assert{c} for each $\alpha$, $\beta\in S^*$, there is a bijection 
$j_{\alpha,\beta}$ on 
$\kappa$ \st\ $j_{\alpha,\beta}\restr r=id_r$, 
$j_{\alpha,\beta}``d_\alpha=d_\beta$, 
$\tilde{j}_{\alpha,\beta}(p_\alpha)=p_\beta$ and 
$\tilde{j}_{\alpha,\beta}(\dotx'_\alpha)=\dotx'_\beta$, where 
$\tilde{j}_{\alpha,\beta}$ is the automorphism on $P$ induced from 
$j_{\alpha,\beta}$ (and the class function on $P$-names induced from this). 
\end{Lemma}
\prf By \Lemmaof{L5}, there are a sequence 
$\seqof{p_\alpha}{\alpha\in S}$ of elements of $P$ below $p$ and a sequence 
$\seqof{\dotx'_\alpha}{\alpha\in S}$ of slim $P$-names \st\ 
$p_\alpha\forces{P}{\dotx_\alpha=\dotx'_\alpha}$ for every $\alpha\in S$. 
Let $d_\alpha=\supp(p_\alpha)\cup\supp(\dotx'_\alpha)$ for $\alpha\in S$. 
Since each $d_\alpha$ is countable and $\kappa$ is $\aleph_0$-inaccessible, 
we can apply Erd\H os-Rado \Thmof{erdoes-rado} to these $d_\alpha$'s and 
obtain a stationary $S'\subseteq S$ \st\ 
$d_\alpha$, $\alpha\in S'$ form a $\Delta$-system with a root $r$. 
We may assume further that 
$d_\alpha\setminus r$, $\alpha\in S'$ have all the same cardinality 
($\leq\aleph_0$). 
For a fixed $\alpha_0\in S'$ and each $\alpha\in S'$ let 
$\mapping{j_\alpha}{\kappa}{\kappa}$ be a bijection \st\ 
$j_\alpha\restr r=id_r$ and $j_\alpha``d_\alpha=d_{\alpha_0}$. 
As before, we 
call the class mapping on $P$-names induced from 
$\tilde{j}_{\alpha}$ also $\tilde{j}_{\alpha}$ for simplicity. 
Now, by $\cardof{Q}<\kappa$ and $\aleph_0$-naccessibility of $\kappa$, we 
have $\cardof{Q}^{\aleph_0}<\kappa$. Hence there are only 
$<\kappa$ many possibilities of $\tilde{j}_{\alpha,\alpha_0}(p_\alpha)$ 
and $\tilde{j}_{\alpha,\alpha_0}(\dotx'_\alpha)$. So we can find a
stationary $S^*\subseteq S'$ \st\ 
$\tilde{j}_{\alpha,\alpha_0}(p_\alpha)$, $\alpha\in S^*$ are all the same 
and also $\tilde{j}_{\alpha,\alpha_0}(\dotx'_\alpha)$, 
$\alpha\in S^*$ are all the same. For $\alpha$, $\beta\in S^*$, it is 
readily seen that $j_{\alpha,\beta}=j_\beta\circ (j_\alpha)^{-1}$ 
satisfies \assertof{c} above. Hence 
$S^*$, 
$\seqof{p_\alpha}{\alpha\in S^*}$ and 
$\seqof{\dotx'_\alpha}{\alpha\in S^*}$ are as required. \qedofLemma 
\qedskip\\ 

\begin{Lemma}\label{L6.5}
Suppose that $P$ is a nice 
side-by-side $\kappa$-product of a \po\ $Q$. 
Then $P$ has the precaliber $(\cardof{Q}^{\aleph_0})^+$. In particular 
$P$ is productively $(\cardof{Q}^{\aleph_0})^+$-cc. 
\end{Lemma}
\prf Let $\lambda=(\cardof{Q}^{\aleph_0})^+$ and suppose that 
$\setof{p_\alpha}{\alpha<\lambda}\subseteq P$. By \Thmof{erdoes-rado}, 
there is $X\in[\lambda]^\lambda$ \st\ $\supp(p_\alpha)$, $\alpha\in X$ form 
a $\Delta$-system with a root $r$. Since 
$\cardof{Q}^{\aleph_0}<\lambda$, we may assume that 
$p_\alpha\restr r$, $\alpha\in X$ are all the same. Hence by 
\Lemmaof{L4.5},\,\assertof{2}, 
$p_\alpha$ and $p_\beta$ are compatible for every $\alpha$, 
$\beta\in X$.\qedofLemma 
\begin{Lemma}\label{L7}
Suppose that $\kappa$ is $\aleph_0$-inaccessible. If $P$ is a proper nice 
side-by-side $\kappa$-product of a \po\ $Q$ with $\cardof{Q}<\kappa$, 
then for any 
sequence $\seqof{p_\alpha}{\alpha\in S}$ of elements of $P$ for a 
stationary $S\subseteq\setof{\alpha<\kappa}{\cof(\alpha)>\omega}$, we have 
$p\forces{P}{\setof{\alpha\in S}{p_\alpha\in\dotG}
	\xmbox{ is stationary subset of }\kappa}$ for some $p\in P$. 
Further, if $\setof{p_\alpha}{\alpha\in S}$ is \st\ 
$\supp(p_\alpha)$, $\alpha\in S$ form a $\Delta$-system with the root 
$r$ and $p_\alpha\restr r$ are all the same for $\alpha\in S$. Then 
$p_\alpha\restr r$ for $\alpha\in S$ can be taken as the $p$ above. 
\end{Lemma}
\prf By the same argument as the proof of \Lemmaof{L6}, it is seen that 
we may suppose that the condition in the last assertion holds: thus 
$\seqof{\supp(p_\alpha)}{\alpha\in S}$ forms a $\Delta$-system with 
the root $r$ and $p_\alpha\restr r=p$ for all $\alpha\in S$. 
We may also assume that $\min(\supp(p_\alpha)\setminus r)\geq\alpha$ for 
every $\alpha\in S$: otherwise, by Fodor's lemma and 
$\cardof{Q}^{\aleph_0}<\kappa$, we can find easily a stationary 
$S^*\subseteq S$ \st\ $p_\alpha$, $\alpha\in S^*$ are all the same. 
Let $q\leq p$ be arbitrary and $\dotC$ be a $P$-name \st\ 
$\forces{P}{\dotC\xmbox{ is a club subset of }\kappa}$. We have to show that 
there is some $q'\leq q$ \st\ 
$q'\forces{P}{\dotC\cap\dotS\not=\emptyset}$ where $\dotS$ is a $P$-name 
for $\setof{\alpha\in S}{p_\alpha\in\dotG}$. Let $\chi$ be sufficiently 
large. 
Since $P$ has the $\kappa$-cc by \Lemmaof{L6.5}, there is an 
$M\prec\calH(\chi)$ \st\ $P$, $Q$, 
$\seqof{p_\alpha}{\alpha\in S},q\in M$;
$M\cap\kappa=\alpha^*\in S$; $\supp(q)\subseteq\alpha^*$; 
$\forces{P}{\min(\dotC\setminus\alpha<\alpha^*)}$ and 
$\supp(p_\alpha)\subseteq\alpha^*$ for every $\alpha<\alpha^*$. Then it 
follows that $\forces{P}{\alpha^*\in\dotC}$. Hence 
$p_{\alpha^*}\forces{P}{\alpha^*\in\dotC\cap\dotS}$. 
We have $\supp(q)\cap\supp(p_{\alpha^*})=r$ and 
$q\restr r\leq p_{\alpha*}\restr r$. Hence $p_{\alpha^*}$ and $q$ are 
compatible by \Lemmaof{L4.5},\,\assertof{2}. 
Let $q'\leq p_{\alpha^*}$, $q$. Then 
$q'\forces{P}{\dotC\cap\dotS\not=\emptyset}$. 
\qedofLemma
\qedskip

Now we are ready to prove \Thmof{main-thm-1}:\medskip\\
\prfof{\Thmof{main-thm-1}}
Let $G$ be a $P$-generic filter over our ground model $V$. To show 
$V[G]\models\HP(\kappa)$, let $f$, $A\in V[G]$ be \st\ 
$\mapping{f}{\kappa}{\psetof{\omega}^{V[G]}}$ and 
$A\subseteq((\psetof{\omega})^{<\aleph_0})^{V[G]}$ is definable. Since 
$P$ is proper, $\omega_1$ in the ground model is preserved. Hence the 
hereditarily countable parameters in a defining formula $\varphi$ for 
$A$ can be 
coded in a subset $a$ of $\omega$ in the way that these parameters can be 
recovered uniformly from $a$ and some ground model sets. 
Hence, by modifying $a$ so that it codes every thing needed, we may 
assume that $\varphi$ contains $a$ as the single parameter. 

Let $\dotf$, $\dota$ and $\dotA$ be 
$P$-names of $f$, $a$ and $A$ respectively \st\ 
$\forces{P}{\mapping{\dotf}{\kappa}{\psetof{\omega}}}$ and 
\[ \forces{P}{\forall x\in(\psetof{\omega})^{<\aleph_0}
	(x\in\dotA\leftrightarrow\calH(\aleph_1)\models\varphi[x,\dota])}. 
\]\noindent
Suppose now that, for $p\in P$, 
\[ p\forces{P}{\mbox{\assertof{h0} in the definition of $\HP(\kappa)$ does 
not hold for }\dotf\mbox{ and }\dotA}. 
\]\noindent
We have to show that $p$ forces the assertion corresponding to 
\assertof{h1} in the definition of 
$\HP(\kappa)$. To do this, let $p'\leq p$ be arbitrary and we show 
that there is $q'\leq p'$ forcing \assertof{h1}. 
By \Lemmaof{L5}, there are $p''\leq p'$ and a slim nice $P$-name 
$\dota'$ \st\ $p''\forces{P}{\dota=\dota'}$. 
By \Lemmaof{L6}, 
there are a 
stationary $S^*\subseteq\setof{\alpha<\kappa}{\cof(\alpha)>\omega}$, a 
sequence $\seqof{\dotx'_\alpha}{\alpha\in S^*}$ of slim nice $P$-names and a 
sequence $\seqof{p_\alpha}{\alpha\in S^*}$ of conditions in $P$ \st\ 
\medskip

\assert{a} $p_\alpha\leq p''$ and 
$p_\alpha\forces{P}{\dotf(\alpha)=\dotx'_\alpha}$ for every $\alpha\in S^*$;

\assert{b} $d_\alpha=\supp(p_\alpha)\cup\supp(\dotx'_\alpha)$,
$\alpha\in S^*$ are all of the same cardinality and form 
a $\Delta$-system with the root $r$. $d_\alpha\cap \supp(\dota')\subseteq r$; 

\assert{c} For each $\alpha$, $\beta\in S^*$ there is a bijection 
$\mapping{j_{\alpha,\beta}}{\kappa}{\kappa}$ \st\ 
$j_{\alpha,\beta}\restr r=id_r$, $j_{\alpha,\beta}``d_\alpha=d_\beta$, 
$\tilde{j}_{\alpha,\beta}(p_\alpha)=p_\beta$ and 
$\tilde{j}_{\alpha,\beta}(\dotx'_\alpha)=\dotx'_\beta$ for every 
$\alpha$, $\beta\in S^*$, where 
$\tilde{j}_{\alpha,\beta}$ is as in \Lemmaof{L6}. 
\medskip \\
Note that, by the last sentence of \assertof{b}, the mapping 
$j_{\alpha,\beta}$ can be chosen 
so that $\tilde j_{\alpha,\beta}(\dota')=\dota'$. 
Let $q=p_\alpha\restr r$ for some/any $\alpha\in S^*$. 
By the inequality in \assertof{a} and since $r$ is the root of 
$d_\alpha$'s, we have $q\leq p''$

Let $\dotS$ be a $P$ name \st\ 
$\forces{P}{\dotS=\setof{\alpha\in S^*}{p_\alpha\in\dotG}}$. By 
\Lemmaof{L7}\medskip

\assert{d} $q \forces{P}{\dotS\xmbox{ is stationary}}$. \medskip\\
Hence, by assumption, we have 
\[ q \forces{P}{\exists n\in\omega(\dotf``(\dotS)^n\not\subseteq\dotA)}. 
\]\noindent
Let $q'\leq q $ pin down the $n$ above as $n^*\in\omega$ and an example 
showing ``$\dotf``(\dotS)^n\not\subseteq\dotA$'' as: 
\begin{markedformula}{$(\ast)$}
\begin{array}[t]{r@{}l}
q'\forces{P}{&\alpha^*_0\tenten\alpha^*_{n^*-1}\in\dotS\ \land\ 
	{\displaystyle\bigwedge_{i<j<n^*}}
		\dotf(\alpha^*_i)\not=\dotf(\alpha^*_j)\\[\jot]
	&\mbox{}\qquad\qquad\land\ \calH(\aleph_1)\models
		\neg\varphi[
		\pairof{\dotf(\alpha^*_0)\tenten\dotf(\alpha^*_{n^*-1})},
			\dota]}
\end{array}
\end{markedformula}
for $\alpha^*_0$\tenten 
$\alpha^*_{n^*-1}\in S^*$. By definition of $\dotS$, 
$p_{\alpha^*_0}\tenten p_{\alpha^*_{n-1}}$ together with $q'$ build a 
sentered set. 
Hence by modifying $q'$, if necessary, we may assume that 
$q'\leq p_{\alpha^*_k}$ for every $k<n^*$. 
Let 
$S^{**}=\setof{\alpha\in S^*}{\supp(q')\cap\supp(p_\alpha)\subseteq r}$. 
Then $S^{**}$ is still stationary subset of $\kappa$. 
For $k<n^*$ and $\alpha\in S^{**}$, let 
$\mapping{j_{k,\alpha}}{\kappa}{\kappa}$ be a bijection \st\ 
$j_{k,\alpha}$ only moves points from 
$(d_{\alpha^*_k}\cup d_\alpha)\setminus r$\,;
$\tilde j_{k,\alpha}(p_{\alpha^*_k})=p_\alpha$\,; and 
$\tilde j_{k,\alpha}(\dotx'_{\alpha^*_k})=\dotx'_\alpha$. Here again 
$\tilde j_{k,\kappa}$ denotes the automorphism on $P$ induced by 
$j_{k,\kappa}$ and the corresponding class mapping on $P$-names. 

Let 
$q^k_\alpha=\tilde j_{k,\alpha}(q')$ for $k<n^*$ and $\alpha\in S^{**}$. 
Then $q^k_\alpha\leq q_\alpha$ for every $\alpha\in S^{**}$. Again 
by \Lemmaof{L7}, 
\[ q\forces{P}{\setof{\alpha\in S^{**}}{q^k_\alpha\in\dotG}
\mbox{ is stationary for every }k<n^*}.
\]\noindent
For $k<n^*$, let $\dotD_k$ be a $P$-name \st\ 
\[ \forces{P}{\dotD_k=\setof{\alpha\in S^*}{q^k_\alpha\in\dotG}}. 
\]\noindent
We claim that 
\[ \begin{array}{r@{}l}
q'\forces{P}{&\calH(\aleph_1)\models\neg\varphi[\pairof{\dotf(\alpha_0)
			\tenten\dotf(\alpha_{n^*-1})},\dota]\\[\jot]
	&\mbox{for every }
		(\alpha_0\tenten\alpha_{n^*-1})\in(\dotD_0\tenten \dotD_{n^*-1})}. 
   
   \end{array}
\]\noindent
Suppose $q''\leq q'$ and $\alpha_0$\tenten $\alpha_{n^*-1}\in\kappa$ be \st\ 
\[q''\forces{P}{(\alpha_0\tenten\alpha_{n^*-1})
	\in(\dotD_0\tenten\dotD_{n^*-1})}.\]\noindent 
Then by definition of 
$\dotD_0$\tenten $\dotD_{n^*-1}$, we have 
$q''\forces{P}{q^0_{\alpha_0}\tenten q^{n^*-1}_{\alpha_{n^*-1}}\in\dotG}$. 
Hence, letting 
$j=j_{0,\alpha_0}\circ\cdots\circ j_{n^*-1,\alpha_{n^*-1}}$ and 
$\tilde j$ the corresponding mapping, 
we have $q''\forces{P}{\tilde j(q')\in\dotG}$. Let $q'''\leq q''$, 
$\tilde{j}(q')$. 
By $(\ast)$, we have 
\[ \tilde j(q')\forces{P}{\calH(\aleph_1)\models
		\neg\varphi[\pairof{\tilde j(\dotf(\alpha^*_0))\tenten
			\tilde j(\dotf(\alpha^*_{n^*-1}))},\tilde j(\dota)]}. 
\]\noindent
Since 
$q''\forces{P}{\tilde j(\dotf(\alpha^*_i))=\tilde j(\dotx'_{\alpha^*_i})
	=\dotx'_{\alpha_i}=\dotf(\alpha_i)}$  
for $i<n^*$ and $\tilde j(\dota)=\dota$, it follows that 
\[ q'''\forces{P}{\calH(\aleph_1)\models
		\neg\varphi[\pairof{\dotf(\alpha_0)\tenten\dotf(\alpha_{n^*-1})},
			\dota]}. 
\]\noindent
\qedof{\Thmof{main-thm-1}}\qedskip
\\
In the proof of \Thmof{main-thm-1}, the definability was not needed in 
the form we introduced it in section 2. What was needed here was merely that 
$A$ can be ``defined'' in some sense from some hereditarily countable 
parameters. E.g., let $\HP^\dagger(\kappa)$ be the principle 
obtained from $\HP(\kappa)$ by replacing the definability by the 
definability by second order property in $(\calH(\aleph_1),\in)$. 
Then \Thmof{main-thm-1} with $\HP^\dagger(\kappa)$ in place of 
$\HP(\kappa)$ still holds. By the argument in the proof of 
\Lemmaof{L1},\,\assertof{4}, this shows that the reals have no second 
order definable well-ordering in the resulting models. C.f.\ 
\cite{abraham-shelah}, where a model of Martin's Axiom is given in which 
there is a $\Delta^2_1$-definable well-ordering of reals. 

Our next goal is the following generalization of \Thmof{main-thm-1}: 
\begin{Thm}\label{main-thm-2} Assume that \CH\ holds. Then, 
for any regular cardinal $\kappa$ and 
any proper nice side-by-side $\kappa$-product of a \po\ $Q$ of 
cardinality $\leq 2^{\aleph_0}$, we have $\forces{P}{\HP(\lambda)}$ for 
every $\aleph_0$-inaccessible $\lambda$. In particular, if 
$\aleph_2\leq\kappa<\aleph_\omega$ then we have 
$\forces{P}{\HP}$. 
\end{Thm}
\begin{Cor}
\mbox{}

\assert{a}
$\MA(countable)$ $+$ $(\clubsuit)$ $+$ \HP\ $+$ $2^{\aleph_0}=\aleph_n$
is consistent relative to \ZFC\ for any $n\in\omega$, $n\geq 1$. 

\assert{b} $\MA(Cohen)$ $+$ \HP\ $+$ $2^{\aleph_0}=\aleph_n$ is 
consistent relative to \ZFC\ for any $n\in\omega$, $n\geq 1$. 
\end{Cor}
\prf 
For $n=1$ the assertions are clear. Suppose $n>1$. 

For \assertof{a}, start from a model of \CH\ and force with 
$\Pi^*_{\aleph_n}\Fn(\omega,2)$. Then the generic extension satisfies 
$2^{\aleph_0}=\aleph_n$, 
$(\clubsuit)$ and 
$\MA(countable)$ by \Lemmaof{544}. By \Thmof{main-thm-2}, this 
model also 
satisfies \HP. 

For \assertof{b}, force with $\Fn(\aleph_n,2)$ over 
a model of \CH. \qedofCor \qedskip

In contrast to Martin's axiom, \HP\ does not imply 
$2^{\aleph_0}=2^{\aleph_1}$: in the combinations of axioms in \Corabove, 
we may also 
require additionally that $2^{\aleph_0}<2^{\aleph_1}$; a model can be 
obtained if we start from a ground model with sufficiently large 
$2^{\aleph_1}$.

Prior to the proof of the theorem, we show the following lemmas which 
should be a sort of folklore. 

We shall call a $P$-name $\dotx$ {\em hereditarily $\lambda$-slim} if 
$\cardof{\setof{\doty,p}{\pairof{\doty,p}\in tcl(\dotx)}}<\lambda$, where 
$tcl(\dotx)$ is the transitive closure of the set $\dotx$. 

\begin{Lemma}\label{L8}
Suppose that $\lambda$ is a regular cardinal, $P$ a $\lambda$-cc \po\ 
and $\chi$ a sufficiently large regular cardinal. 
If $M\prec\calH(\chi)$ is \st\ 
$P,\lambda+1\in M$ and $[M]^{<\lambda}\subseteq M$, then:

\assert{1} $P\cap M$ is a regular subordering of $P$.\smallskip

\assert{2} For any $P$-name $\dotx$, there is a 
hereditarily 
$\lambda$-slim $P$-name $\doty$ \st\ 
$\forces{P}{\dotx\in\calH(\lambda)\rightarrow\dotx=\doty}$. 
\smallskip

\assert{3} If $\dotx\in M$ is a hereditarily $\lambda$-slim $P$-name, 
then $\dotx$ is a $P\cap M$-name. \smallskip

\assert{4} If $\dotx$ is a hereditarily $\lambda$-slim $P\cap M$-name 
then $x\in M$.\smallskip

\assert{5} If $p\in P\cap M$, $\varphi$ is a first-order formula in the 
language of 
set-theory and $\dotx_0$\tenten $\dotx_{n-1}$ are hereditarily 
$\lambda$-slim $P$-names in M \st\ 
\[ \forces{P}{\dotx_0\tenten\dotx_{n-1}\in\calH(\lambda)}
\]\noindent
then we have 
\begin{markedformula}{$(*)$}
p\forces{P}{\calH(\lambda)\models\varphi[\dotx_0\tenten\dotx_{n-1}]}
	 \equivto
	p\forces{P\cap M}{%
		\calH(\lambda)\models\varphi[\dotx_0\tenten\dotx_{n-1}]}.
\end{markedformula}

\assert{6} If $G$ is a $P$-generic filter and $G_M=G\cap M$ then $G_M$ is 
a $P\cap M$-generic filter and 
\[ V[G]\models\calH(\lambda)^{V[G_M]}\prec\calH(\lambda). 
\]\noindent
\end{Lemma}
\prf
\assertof{1}\,: Let $A\subseteq P\cap M$ be a maximal antichain. By the 
$\lambda$-cc of $P$, we have $\cardof{A}<\lambda$. By assumption on 
$M$ it follows that $A\in M$. Since 
$M\models\mbox{``}A\xmbox{ is a maximal antichain in }P\mbox{''}$, 
$A$ is really a maximal antichain in $P$ by elementarity. 

\assertof{2} We show this by induction on $rnk(\dotx)$. Let $\dotx$ be a 
$P$-name and assume that we have 
shown the assertion for every $P$-name with rank lower than that of 
$\dotx$. Let $\dotx=\setof{\pairof{\dotx_\alpha, p_\alpha}}{\alpha<\kappa}$ for 
some cardinal $\kappa$. By assumption there are hereditarily 
$\lambda$-slim $P$-names $\dotx'_\alpha$ for $\alpha<\kappa$ \st\ 
$\forces{P}{\dotx_\alpha\in\calH(\lambda)
	\rightarrow\dotx_\alpha=\dotx'_\alpha}$. 
By the $\lambda$-cc of $P$, there are $\mu<\lambda$ and a $P$-name 
$\dotf$ \st\ 
\[\forces{P}{\mapping{\dotf}{\mu}{\dotx}\land(\dotx\in\calH(\lambda)
	\rightarrow\dotf\mbox{ is onto})}.\]\noindent
For each $\beta<\mu$, let $D_\beta\subseteq P$ be a maximal antichain 
\st\ for every $r\in D_\beta$ there is $\alpha_{r,\beta}<\kappa$ \st\ 
$r\forces{P}{p_{\alpha_{r,\beta}}\in\dotG}$ and 
$r\forces{P}{\dotf(\beta)=\dotx'_{\alpha_{r,\beta}}}$. By the 
$\lambda$-c.c.\ of $P$, we have $\cardof{D_\beta}<\lambda$. Let 
\[ \doty=\setof{\pairof{\dotx'_{\alpha_{r,\beta}},p_{\alpha_{r,\beta}}}}% 
	{r\in D_\beta\mbox{ for some }\beta<\mu}. 
\]\noindent
Then $\doty$ is as desired. 

\assertof{3}\,: Suppose that $\dotx\in M$ is a hereditarily 
$\lambda$-slim $P$-name. Then 
$U=\setof{p\in P}{\pairof{\doty, p}\in tcl(\dotx)}$ is an element of $M$ and of 
cardinality $<\lambda$. Hence $U\subseteq P\cap M$ by assumption 
on $M$. It follows that $\dotx$ is a $P\cap M$-name. 

\assertof{4}\,: By induction on $rnk(\dotx)$. Suppose that $\dotx$ is a 
hereditarily $\lambda$-slim $P\cap M$-name. Say 
$\dotx=\setof{\pairof{\dotx_\alpha,p_\alpha}}{\alpha<\mu}$ for some 
$\mu<\lambda$. By induction hypothesis and since each $\dotx_\alpha$ is 
also hereditarily $\lambda$-slim, we have $\dotx\subseteq M$. Hence by 
$\cardof{\dotx}\leq\mu<\lambda$ and by assumption on $M$ it follows that 
$\dotx\in M$. 

\assertof{5}\,: We show this by induction on $\varphi$. 
If $\varphi$ is atomic, then the equivalence $(*)$ clearly holds. The 
induction steps for $\neg$, $\land$, $\lor$ are also clear. Assume that 
$\varphi=\exists x\psi(x,x_0\tenten x_{n-1})$
and $(*)$ holds for 
$\psi$. Suppose first that 
\[p\forces{P}{\calH(\lambda)\models\varphi[\dotx_0\tenten\dotx_{n-1}]}. 
\]
Since $p$, $\lambda$, 
$\dotx_0$\tenten $\dotx_{n-1}\in M$, 
\[ M\models p\forces{P}{\calH(\lambda)
		\models\varphi[\dotx_0\tenten\dotx_{n-1}]}
\]\noindent
by elementarity of $M$. Hence there is a $P$-name $\dotx\in M$ \st\ 
\[M\models p\forces{P}{\calH(\lambda)
		\models\psi[\dotx,\dotx_0\tenten\dotx_{n-1}]}. 
\]\noindent
Hence 
\[p\forces{P}{\calH(\lambda)
		\models\psi[\dotx,\dotx_0\tenten\dotx_{n-1}]}. 
\]\noindent
By \assertof{2}, we may assume that $\dotx$ is a hereditarily $\lambda$-slim 
$P$-name. Then, by \assertof{3}, $\dotx$ is a $P\cap M$-name. By 
induction hypothesis, it follows that 
\[p\forces{P\cap M}{\calH(\lambda)
		\models\psi[\dotx,\dotx_0\tenten\dotx_{n-1}]}
\]\noindent
and hence 
$p\forces{P\cap M}{\calH(\lambda)
	\models\varphi[\dotx_0\tenten\dotx_{n-1}]}$. 

Suppose now that 
\[p\forces{P\cap M}{\calH(\chi)\models\varphi[\dotx_0\tenten \dotx_{n-1}]}. 
\]\noindent
Then, by \assertof{2}, there is a hereditarily $\lambda$-slim 
$P\cap M$-name $\dotx$ \st\ 
\[p\forces{P\cap M}{%
	\calH(\lambda)\models\psi[\dotx,\dotx_0\tenten\dotx_{n-1}]}. \]
By \assertof{4}, $\dotx\in M$. Hence by induction hypothesis 
\[p\forces{P}{%
	\calH(\lambda)\models\psi[\dotx,\dotx_0\tenten\dotx_{n-1}]}. \]
Thus we have 
$p\forces{P}{\calH(\lambda)\models\varphi[\dotx_0\tenten\dotx_{n-1}]}$. 

\assertof{6}\,: $G_M$ is a $P\cap M$-generic filter by \assertof{1}. The 
rest of the assertion follows immediately from \assertof{5}. 
\qedofLemma\qedskip

For nice side-by-side products, \assertof{5} and \assertof{6} of the 
lemma above 
can be still slightly improved. If $P$ is a nice side-by-side 
$\kappa$-product of some $Q$ and $\dotx$ is a $P$-name, let 
$\supp(\dotx)=
	\bigcup\setof{\supp(p)}{\pairof{\doty,p}\in tcl(\dotx)\xmbox{ for some }\doty}$. 
Note that, if $\supp(\dotx)\subseteq X$ for $X\subseteq\kappa$, then 
$\dotx$ is a $P_X$ name where $P_X$ is defined as in \assertof{n4} in the 
definition of nice side-by-side product. 
\begin{Lemma}\label{L9}
Suppose that $P$ is a nice side-by-side $\kappa$-product of a \po\ $Q$ 
with $Q_\alpha$ and $\iota_\alpha$ for $\alpha<\kappa$ as in the 
definition of nice side-by-side product. For a regular 
cardinal $\lambda>\cardof{Q}^{\aleph_0}$ and 
$X\in [\kappa]^{\lambda}$, let $P_X$ be the regular subordering of 
$P$ as defined in \assertof{n3}. Then: 
\smallskip

\assert{1} For any $p\in P_X$, first-order formula $\varphi$ in the 
language of 
set-theory and hereditarily 
$\lambda$-slim $P_X$-names $\dotx_0$\tenten $\dotx_{n-1}$ 
\begin{markedformula}{$(**)$}
p\forces{P}{\calH(\lambda)\models\varphi[\dotx_0\tenten\dotx_{n-1}]}
	 \equivto
	p\forces{P_X}{\calH(\lambda)\models\varphi[\dotx_0\tenten\dotx_{n-1}]}.
\end{markedformula}
\par
\assert{2} If $G$ is a $P$-generic filter and $G_X=G\cap P_X$ then 
\[ V[G]\models\calH(\lambda)^{V[G_X]}\prec\calH(\lambda). 
\]\noindent
\end{Lemma}
\prf
\assertof{1}\,: By induction on $\varphi$. If $\varphi$ is atomic, $(**)$ 
clearly holds. The induction steps for $\neg$, $\land$, $\lor$ are easy. 
Assume that $\varphi$ is $\exists x\psi(x,x_0\tenten x_{n-1})$ and we 
have proved $(**)$ for $\psi$. 
First, suppose that we have 
\[p\forces{P}{\calH(\lambda)\models\varphi[\dotx_0\tenten\dotx_{n-1}]}. 
\]\noindent
$P$ has the $(\cardof{Q}^{\aleph_0})^+$-cc by \Lemmaof{L6.5}. 
Hence, by \Lemmaof{L8},\,\assertof{2}, there is a hereditarily 
$\lambda$-slim $P$-name $\dotx$ \st\ 
\[p\forces{P}{\calH(\lambda)\models\psi[\dotx, \dotx_0\tenten\dotx_{n-1}]}. 
\]\noindent
Since $\cardof{\supp(\dotx)}<\lambda$, there is a bijection 
$\mapping{j}{\kappa}{\kappa}$ \st\ $j``\supp(\dotx)\subseteq X$, 
$\tilde j(p)=p$ and 
$\tilde{j}(\dotx_i)=\dotx_i$ for all $i<n$ where $\tilde j$ is as in 
\assertof{n4}. Then $\tilde{j}(\dotx)$ is a 
$P_X$-name and
\[p\forces{P}{\calH(\lambda)\models
		\psi[\tilde{j}(\dotx),\dotx_0\tenten\dotx_{n-1}]}. 
\]\noindent
By induction hypothesis, it follows that 
\[ p\forces{P_X}{\calH(\lambda)\models
		\psi[\tilde{j}(\dotx),\dotx_0\tenten\dotx_{n-1}]}
\]\noindent
and hence
$p\forces{P_X}{\calH(\lambda)\models
		\exists x\psi[x,\dotx_0\tenten\dotx_{n-1}]}$. 

Suppose now that 
\[ p\forces{P_X}{\calH(\lambda)\models
		\exists x\psi[x,\dotx_0\tenten\dotx_{n-1}]}. 
\]\noindent
Then there is a $P_X$-name $\dotx$ \st\
\[ p\forces{P_X}{\calH(\lambda)\models
		\psi[\dotx,\dotx_0\tenten\dotx_{n-1}]}. 
\]\noindent
By induction hypothesis, it follows that 
\[ p\forces{P}{\calH(\lambda)\models
		\psi[\dotx,\dotx_0\tenten\dotx_{n-1}]}
\]\noindent
and hence 
$p\forces{P}{\calH(\lambda)\models
		\exists x\psi[x,\dotx_0\tenten\dotx_{n-1}]}$.
\smallskip

\assertof{2} follows immediately from \assertof{1}. 
\qedofLemma
\qedskip\\
\prfof{\Thmof{main-thm-2}}
Let $\kappa$ be regular 
$P$ be a proper nice side-by-side $\kappa$-product of a partial ordering 
$Q$ of cardinality $\leq 2^{\aleph_0}=\aleph_1$. By \Lemmaof{L6.5}, 
$P$ has the $\aleph_2$-cc. By this and properness, $P$ preserves 
cofinality and cardinals. Since $\cardof{P}\leq\kappa$ we have 
$2^{\aleph_0}\leq\kappa$. Hence by \Lemmaof{L1},\,\assertof{3}, we have 
$\forces{P}{\HP(\lambda)}$ for every regular $\lambda>\kappa$. Let 
$\aleph_2\leq\lambda\leq\kappa$ be $\aleph_0$-inaccessible 
and suppose that $p\in P$ and, $\dotf$ and $\dotx$ are 
$P$-names \st\ 
\[ p\forces{P}{\mapping{\dotf}{\lambda}{\psetof{\omega}},\,
	\dotx\in\calH(\aleph_1)}
\]\noindent
We would like to show:
\[ p\forces{P}{\mbox{\assertof{h0} or \assertof{h1} in the definition of }
		\HP(\lambda)\mbox{ holds for }\dotf\mbox{ and }\dotA}
\]\noindent
where $\dotA$ is a $P$-name \st\ for a first order formula $\varphi$ in 
the language of set theory and 
\[ p\forces{P}{\dotA
	=\setof{u\in\psetof{\omega}}{\calH(\aleph_1)\models\varphi[u,\dotx]}}.
\]\noindent
By \Lemmaof{L8},\,\assertof{2}, we may assume that $\dotf$ and 
$\dotx$ are $\lambda$-slim names. Hence there is 
$X\subseteq[\kappa]^\lambda$ \st\ $p\in P_X$ and, $\dotf$ and $\dotx$ are 
$P_X$-names. 
Since $P_X$ is a proper nice side-by-side $\lambda$-product of $Q$ and 
$\lambda$ is $\aleph_0$-inaccessible by assumption, 
we have $\forces{P_X}{\HP(\lambda)}$ by \Thmof{main-thm-1}. Hence 
\[ 
p\forces{P_X}{\mbox{\assertof{h0} or \assertof{h1} in the definition of }
		\HP(\lambda)\mbox{ holds for }\dotf\mbox{ and }\dotA'}
\]\noindent
where $\dotA'$ is a $P_X$-name \st\ 
\[\forces{P_X}{\dotA'
	=\setof{u\in\psetof{\omega}}{\calH(\aleph_1)\models\varphi[u,\dotx]}}.
\]\noindent
By \Lemmaof{L9} we have 
$p\forces{P}{\dotA'=\dotA\cap\calH(\aleph_1)^{V[G_X]}}$. 
Since 
$\forces{P_X}{P_{\kappa\setminus X}\xmbox{ has the }\lambda\xmbox{-c.c.}}$, 
by \Lemmaof{L6.5}, 
stationary subsets of $\lambda$ in the generic extension by $P_X$ are 
still stationary in the generic extension by $P$. Hence it follows that 
\[\forces{P}{\mbox{\assertof{h0} or \assertof{h1} in the definition of }
		\HP(\lambda)\mbox{ holds for }\dotf\mbox{ and }\dotA}. 
\]\noindent
This shows $\forces{P}{\HP(\lambda)}$. 
\qedof{\Thmof{main-thm-2}}\qedskip

Now, let us turn to the models of \IP. We need a well-known fact about 
product of \pos\ which is formulated here in the following lemma for nice 
side-by-side products:
\begin{Lemma}\label{L10}
Let $P$ be a nice side-by-side $\kappa$-product of some \po\ $Q$. 
Suppose that $X,Y\subseteq\kappa$ and $Z=X\cap Y$. For a $P$-generic 
filter $G$ over $V$, let $G_X=G\cap P_X$. By 
\Lemmaof{L4.5},\,\assertof{0} $G_X$ is a $P_X$-generic filter over $V$. 
Let $G_Y$ and $G_Z$ be defined similarly. 
If $x\in V[G_X]\cap\On$ is not an element of 
$V[G_Z]$ then $x$ is neither an element of $V[G_Y]$. 
\end{Lemma}
\prf

By replacing $V$ with $V[G_Z]$ and $P$ with 
$P_{\kappa\setminus Z}$ we may assume that $Z=\emptyset$. Let 
$y\in V[G_Y]\cap\On$. We show that $x\not=y$. Let $\dotx$ be a $P_X$ name 
of $x$ and $\doty$ a $P_Y$-name of $y$. For arbitrary $p\in G$, we have 
$p\restr X\notforces{P_X}{\dotx\in V}$. Hence there is $\alpha\in\On$ 
\st\ $p\restr X$ does not decide ``$\alpha\in\dotx$'' in $P_X$. Let 
$q\in P_Y$ be \st\ $q\leq p\restr Y$ and $q$ decide 
``$\alpha\in\doty$'', say $q\forces{P_Y}{\alpha\in\doty}$. Let 
$p'\in P_X$ be \st\ $p'\leq p\restr X$ and 
$p'\forces{P_X}{\alpha\not\in\dotx}$. Then $p'\land q\leq p$ and 
$p'\land q\forces{P}{\dotx\not=\doty}$. \qed
\begin{Thm}\label{main-thm-3}
Let $\kappa$ be an $\aleph_0$-inaccessible cardinal and 
$\lambda\leq\kappa$. Suppose that $P$ is a proper nice side-by-side 
$\kappa$-product of a \po\ $Q$ \st\ $\cardof{Q}^{\aleph_0}<\kappa$. 
If $P$ has the $\lambda$-cc, then $\forces{P}{\IP(\kappa,\lambda)}$. 
\end{Thm}
\prf Let $G$ be a $P$-generic filter over $V$. In $V[G]$, let 
$\mapping{f}{\kappa}{\psetof{\omega}}$ and  
$\mapping{g}{(\psetof{\omega})^{<\aleph_0}}{\psetof{\omega}}$ be definable, 
say by a formula $\varphi$. As in the proof of \Thmof{main-thm-1}, we may 
assume that 
$\varphi$ has a real $a\in V[G]$ as  its unique parameter. 
Let $\dotf$, $\dota$ and $\dotg$ be $P$-names of $f$, $a$ and $g$ \st\ 
$\forces{P}{\mapping{\dotf}{\kappa}{\psetof{\omega}}}$, 
$\forces{P}{\mapping{\dotg}{(\psetof{\omega})^{<\aleph_0}}{\psetof{\omega}}}$ 
and 
\[\forces{P}{\forall\overline{x}\in(\psetof{\omega})^{<\aleph_0}
	\forall x\in\psetof{\omega}\bigl({\dotg(\overline{x})=x\leftrightarrow
	\calH(\aleph_1)\models\varphi(\overline{x},x,\dota)}\bigr)}. 
\]\noindent
Suppose that, for a $p\in G$, we have 
\[ p\forces{P}{\mbox{\assertof{i0} does not hold for }\dotf\mbox{ and }\dotg}.
\]\noindent
In particular we have
\[ p\forces{P}{%
	\forall \alpha<\kappa(\setof{\beta\in\kappa}{\dotf(\beta)=\dotf(\alpha)}
	\mbox{ is non-stationary})}.
\]\noindent
\begin{Claim}
There is a stationary $S\subseteq\kappa$ \st\ 
\[ p\forces{P}{\dotf\restr S\mbox{ is 1-1}}.
\]\noindent
\end{Claim}
\prfofClaim
By assumption and the $\lambda$-cc of $P$ there are club sets 
$C_\alpha\subseteq\kappa$ for $\alpha<\kappa$ \st\
\[ p\forces{P}{%
	C_\alpha\cap\setof{\beta\in\kappa}{\dotf(\beta)=\dotf(\alpha)}=\emptyset}.
\]\noindent
Then $S=\Delta_{\alpha}C_\alpha$ is club and hence stationary, and has the 
desired property.\qedofClaim\qedskip\\
We show that $p$ forces \assertof{i1} for these $\dotf$ and 
$\dotg$. Let $p'\leq p$ be arbitrary. It is enough to show that there is 
$p^*\leq p'$ forcing \assertof{i1}. 
Similarly to the proof of \Thmof{main-thm-1}, 
we can find a $p''\leq p'$, a slim nice $P$-name $\dota'$ of a real, 
a stationary $S^*\subseteq S$, a sequence 
$\seqof{\dotx'_\alpha}{\alpha\in S^*}$ of slim nice $P$-names and a 
sequence $\seqof{p_\alpha}{\alpha\in S^*}$ of conditions in $P$ \st\ 
\medskip

\assert{a} $p''\forces{P}{\dota=\dota'}$
$p_\alpha\leq p''$ and 
$p_\alpha\forces{P}{\dotf(\alpha)=\dotx'_\alpha}$ for every $\alpha\in S^*$;

\assert{b} $d_\alpha=\supp(p_\alpha)\cup\supp(\dotx'_\alpha)$,
$\alpha\in S^*$ are all of the same cardinality and form 
a $\Delta$-system with a root $r$. 
$\supp(\dota')\cap d_\alpha\subseteq r$ for every $\alpha\in S^*$; 

\assert{c} for each $\alpha$, $\beta\in S^*$ there is a bijection 
$\mapping{j_{\alpha,\beta}}{\kappa}{\kappa}$ \st\ 
$j_{\alpha,\beta}\restr\kappa\setminus(d_\alpha\Delta d_\beta)$ is an 
identity mapping, $j_{\alpha,\beta}``d_\alpha=d_\beta$, 
$\tilde{j}_{\alpha,\beta}(p_\alpha)=p_\beta$ and 
$\tilde{j}_{\alpha,\beta}(\dotx'_\alpha)=\dotx'_\beta$ for every 
$\alpha$, $\beta\in S^*$, where 
$\tilde{j}_{\alpha,\beta}$ is as in \Lemmaof{L6};\medskip
\\
By \assertof{c} $p_\alpha\restr r$ for $\alpha\in S^*$ are all the same 
and is stronger than $p''$. 
Let $q=p_\alpha\restr r$ for some/any $\alpha\in S^*$ and $\dotS$ be a 
$P$ name \st\ 
$\forces{P}{\dotS=\setof{\alpha\in S^*}{p_\alpha\in\dotG}}$. 
By \Lemmaof{L7}, it follows that\medskip

\assert{d} $q \forces{P}{\dotS
	\xmbox{ is stationary}}$. \medskip\\
Hence, by assumption, 
\[ q \forces{P}{%
	\exists n\in\omega(\cardof{\dotg``\dotf``(\dotS)^n}\geq\lambda)}. 
\]\noindent
Let $q'\leq q $ pin down the $n$ above as $n^*\in\omega$. 
Let $S^{**}=\setof{\alpha\in S^*}{\supp(q')\cap\supp(p_\alpha)\subseteq r}$. 
Then $S^{**}$ is still stationary subset of $\kappa$. 
Let $\pairof{\alpha_0\tenten\alpha_{n^*-1}}\in (S^{**})^{n^*}$ and 
$q''=q\land p_{\alpha_0}\land\cdots\land p_{\alpha^*-1}$. 
\begin{Claim}
\[ q''\notforces{P}{%
	\dotg(\dotf(\alpha_0)\tenten\dotf(\alpha_{n^*-1}))\in V[\dotG_r]} 
\]\noindent
\end{Claim}
\prfofClaim
Otherwise,   
by the $\lambda$-cc of $P$, there is a $P_r$-name $\dotX$ \st\ 
$\forces{P_r}{\cardof{\dotX}<\lambda}$ and 
\[q''\forces{P}{\dotg(\dotf(\alpha_0)\tenten\dotf(\alpha_{n^*-1}))\in\dotX}.
\]\noindent
Hence, by \assertof{c}, it follows that 
$q''\forces{P}{\dotg``\dotf``(\dotS)^n\subseteq \dotX}$. This is a 
contradiction to the choice of $n^*$ and $q''$. 
\qedofClaim\qedskip

Let $p^*\leq q''$ be \st\ 
\[p^*\forces{P}{\dotg(\dotf(\alpha_0)\tenten\dotf(\alpha_{n^*-1}))
	\not\in V[\dotG_r]}.\]\noindent
By shrinking $S^{**}$, if necessary, we may assume that 
$\supp(p^*)\cap\supp(p_\alpha)\subseteq r$ for all $\alpha\in S^{**}$. 
%% Let 
%% $X=\supp(p^*)\setminus\bigcup_{i<n^*}\supp(p_{\alpha_i})$. 
For 
$i<n^*$, let $\dotD_i$ be a $P$-name \st\ 
%\[ \forces{P}{\dotD_i=\setof{\alpha\in S^{**}}{%
%	\tilde{j}_{\alpha_i,\alpha}(p^*\restr (X\cup\supp(p_{\alpha_i})))
%		\in\dotG }}.
%\]\noindent
\[ \forces{P}{\dotD_i=\setof{\alpha\in S^{**}}{%
	\tilde{j}_{\alpha_i,\alpha}(p^*)
		\in\dotG }}.
\]\noindent
Then, by \Lemmaof{L7}, we have 
$p^*\forces{P}{\dotD_i\xmbox{ is a stationary subset of }\kappa}$ for all 
$i<n^*$. By definition of $\dotD_i$'s, we have
\[ 
\begin{array}{r@{}l}
p^*\forces{P}{&\forall\beta_0\cdots\forall\beta_{n^*-1}
	(\pairof{\beta_0\tenten\beta_{n^*-1}}\in(\dotD_0\tenten \dotD_{n^*-1})
			\\[\jot]
	&\,\rightarrow
		\dotg(\dotf(\beta_0)\tenten\dotf(\beta_{n^*-1}))
			\not\in V[\dotG_r])}. 
\end{array}
\]\noindent
In general, the support of $P$-names $\doty(\beta_0\tenten\beta_{n^*-1})$ of  
$\dotg(\dotf(\beta_0)\tenten\dotf(\beta_{n^*-1}))$ for 
$\pairof{\beta_0\tenten\beta_{n^*-1}}\in(S^{**})^{n^*}$ can be anything.  But 
since $g$ is forced to be definable by $\varphi$, we can apply 
\Lemmaof{L9} to choose these names so that 
\[\supp(\doty(\beta_0\tenten\beta_{n^*-1}))
\cap\supp(\doty(\beta'_0\tenten\beta'_{n^*-1}))\subseteq r,\]
\noindent 
whenever $\beta_0$\tenten $\beta_{n^*-1}$, $\beta'_0$\tenten 
$\beta'_{n^*-1}\in S^{**}$ are pairwise distinct. 
Hence, by \Lemmaof{L10}, we conclude that 
$p^*\forces{P}{\mbox{\assertof{i1} holds.}}$
\qedof{\Thmof{main-thm-3}}\qedskip

Similarly to the proof of \Thmof{main-thm-2}, we can now proof 
the following: 
\begin{Thm}
Suppose that $\kappa$ is a regular cardinal and $P$ is a 
proper nice side-by-side $\kappa$-product of a \po\ $Q$ of 
cardinality $\leq 2^{\aleph_0}$. Then for any 
$\aleph_0$-inaccessible $\lambda$ we have 
$\forces{P}{\IP(\lambda,\aleph_2)}$. 
In particular, if \GCH\ holds and $\aleph_2\geq\kappa<\aleph_\omega$, 
then $\forces{P}{\IP}$. If further $P$ satisfies the ccc then 
$\forces{P}{\IP^+}$. \qed
\end{Thm}

\begin{Cor}\mbox{}

\assert{a}
$\MA(countable)$ $+$ $(\clubsuit)$ $+$ \HP\ $+$ \IP\ $+$ 
$2^{\aleph_0}=\aleph_n$ is consistent relative 
to \ZFC\ for every $n\in\omega$, $n\geq 1$. \smallskip

\assert{b}
$\MA(Cohen)$ $+$ \HP\ $+$ $\IP^+$ $+$ $2^{\aleph_0}=\aleph_n$ is 
consistent relative to \ZFC\ for every $n\in\omega$, $n\geq 1$. \qed 
\end{Cor}

Finally, we would like to mention the following open problems: 
\begin{Problem}
Does  $\IP(\kappa)$ imply $\HP(\kappa)$ or vice versa~? 
\end{Problem}
\begin{Problem}
Is $\HP(\aleph_{\omega+1})$ (or $\IP(\aleph_{\omega+1})$) consistent with 
$2^{\aleph_0}>\aleph_\omega$~?
\end{Problem}
Note that the second problem could be also formulated as whether 
$\HP$ (or $\IP$) implies that $2^{\aleph_0}<\aleph_\omega$. 
\bigskip\\
{\bf Acknowledgment.} I would like to thank Lajos Soukup and Yoshihiro 
Abe for pointing out several minor errors in an earlier version of the paper. 

\newpage

%--------------------------

\mbox{}\vfill\mbox{}\hfill
\parbox[t]{7cm}{
\noindent
{\bf Author's address}\bigskip\bigskip\\
\it
{\rm Saka\'e Fuchino}\smallskip\\
Institut f\"ur Mathematik II,\\ Freie Universit\"at Berlin\\
14195 Berlin, Germany\medskip\\
{\tt fuchino@math.fu-berlin.de}}
\end{document}